\documentclass{amsart}
\usepackage{kotex, amsmath, amsfonts, amssymb, amsthm, epstopdf, latexsym, enumerate, mathrsfs, mathtools, ytableau}
\usepackage{wrapfig}

\newtheorem{Thm}{Theorem}
\newtheorem{Pro}[Thm]{Proposition}
\newtheorem{Lem}[Thm]{Lemma}
\newtheorem{Cor}[Thm]{Corollary}

\newtheorem{Con}[Thm]{Conjecture}

\theoremstyle{plain}

\newtheorem*{Cla*}{Claim}

\theoremstyle{definition}
\newtheorem{Def}[Thm]{Definition}

\theoremstyle{remark}
\newtheorem*{rem}{Remark}

\newcommand{\mb}[1]{\mathbb{#1}}

\title{Explicit formulas for $e$-positivity of chromatic quasisymmetric functions}

\author{Seung Jin Lee, Sue Kyong Y. Soh}

\begin{document}
\begin{abstract}
	\par
	In 1993, Stanley and Stembridge conjectured that a chromatic symmetric function of any $(3+1)$-free poset is $e$-positive. Guay-Paquet reduced the conjecture to $(3+1)$- and $(2+2)$-free posets which are also called natural unit interval orders. Shareshian and Wachs defined chromatic quasisymmetric functions, generalizing chromatic symmetric functions, and conjectured that a chromatic quasisymmetric function of any natural unit interval order is $e$-positive and $e$-unimodal. \\

For a given natural interval order, there is a corresponding partition $\lambda$ and we denote the chromatic quasisymmetric function by $X_\lambda$. The first author introduced local linear relations for chromatic quasisymmetric functions. In this paper, we prove a powerful generalization of the above-mentioned local linear relations, called a rectangular lemma, which also generalizes the formula in \cite{huh2020melting}. Such a lemma can be applied to describe explicit formulas for $e$-positivity of a chromatic symmetric function $X_\lambda$ where $\lambda$ is contained in a rectangle. We also suggest some conjectural formulas for $e$-positivity when $\lambda$ is not contained in a rectangle by applying the rectangular lemma.
\end{abstract}
\maketitle

\section{Introduction}

    $\quad$	In 1995, Stanley \cite{stanley1995symmetric} introduced a chromatic symmetric function $X_G(x)$, which is a generalization of a chromatic polynomial of $G$, associated to a simple graph $G$. There has been plenty of research about the chromatic symmetric function in diverse areas. Stanley and Stembridge \cite{stanley1993immanants} introduced one of the famous conjectures on chromatic symmetric functions; a chromatic symmetric function of any $(3+1)$-free poset is a linear sum of elementary symmetric functions $\{e_{{\lambda}}\}$-basis with non-negative coefficients. Guay-Paquet \cite{guay2013modular} in 2014 proved that if a chromatic symmetric function of a natural unit interval order set is $e$-positive, then Stanley and Stembridge's conjecture is also true. Shareshian and Wachs \cite{shareshian2016chromatic} introduced a chromatic quasisymmetric refinement of Stanley’s chromatic symmetric function in 2016. They conjectured a chromatic quasisymmetric function of any natural unit interval order is $e$-positive and $e$-unimodal. That is, this conjecture is a refinement of Stanley and Stembridge's conjecture.\\
	
	In 2018, the first author \cite{lee2018linear} introduced local linear relations on unicellular LLT polynomials, and the chromatic quasisymmetric functions. In 2020, Huh, Nam, and Yoo \cite{huh2020melting} generalized and utilized the linear relations of the chromatic quasisymmetric functions from \cite{lee2018linear} to find expanded local linear relations on the chromatic quasisymmetric functions. In this paper, we prove a rectangular lemma, which is a vast generalization of above-mentioned relations. We can partially apply the rectangular lemma although $\lambda$ does not contained in a rectangle. If $\lambda$ is contained in a rectangle of size $d_1\times d_2$ satisfying $d_2\geq d_1$, then the rectangular lemma suggests that $X_\lambda$ can be written explicitly as a positive linear combination in terms of $X_{d_2^r}$ where $0\leq r \leq d_1$. For example, when $\lambda=(2,2)$, $d_1=2$, and $d_2=n-2$, the rectangular lemma says 

\[\text{\scalebox{0.83}{$X_\lambda=\frac{1}{[n-2]_q[n-3]_q}\Big([n-4]_q[n-5]_q X_{(0)}+q^{n-5}[2]_q[2]_q[n-4]_q X_{(n-2)}+q^{2n-8}[2]_q X_{(n-2,n-2)}\Big)$}}\] where $n\geq4$.

We also provide explicit formulas for $e$-positivity of the chromatic quasisymmetric functions when $\lambda$ is contained in a rectangle. This $e$-positivity was described in \cite{abrue2020modularlaw,cho2019positivity, harada2019cohomology} previously, but the fact that the rectangular lemma can be applied to a broader class of partitions is useful to study their $e$-positivity of chromatic quasisymmetric functions. We will prove some results and make conjectures for some cases. Moreover, when $\lambda$ is contained in a rectangle, one can compare our explicit formula with the formula in terms of $q$-hit numbers and rook placements, proved by Abreu and Nigro \cite{abrue2020modularlaw}. Then we show that our formula also provides an explicit formula for $q$-hit numbers (See section 4). \\

Also, combinatorial objects and quantities appearing in the rectangular lemma also appeared in different papers \cite{Jangsoo2014,Shi2012} studying maximal parabolic Kazhdan-Lusztig polynomials and Dyck tiling. It would be interesting if the coefficients in $e$-expansion of $X_\lambda$ is related to quantities appearing in their work.\\

	The contents of the paper are organized as follows. Section 2 lays out necessary definitions, notations, and known results used in Sections 3 and 4. In Section 3, we introduce the rectangular lemma. Moreover, we show that a similar formula holds for $e$-positivity of chromatic quasisymmetric functions when $\lambda$ is contained in a rectangle. In Section 4, we discuss the relationship between our results and $q$-hit numbers, rook placements, and other results. In Section 5, we apply our rectangular lemma to find a few explicit formulas for $e$-positivity of chromatic quasisymmetric functions when $\lambda$ is a partition such that $\lambda$ is contained in a rectangle except the first row of $\lambda$. We also suggest some conjectural formulas for $e$-positivity of chromatic quasisymmetric functions when $\lambda$ is not abelian.

\section{Preliminaries}
In this section, we state definitions, notations, and known results needed in this paper. More details can be found in \cite{alexandersson2018llt, carlsson2018proof, haglund2005combinatorial, huh2020melting, lascoux1997ribbon, lee2018linear, shareshian2016chromatic, stanley1995symmetric}.\\

For a positive integer $n$, we set $[n]=\{1,2,\cdots,n\}$ and $q$-integer \[[n]_q=1+q+q^2+\cdots+q^{n-1}.\]

\subsection{Chromatic quasisymmetric functions}

A simple graph is a graph with no loops and no multiple edges. A \emph{proper coloring} of a simple graph $G=({V},{E})$ with a vertex set ${V}$ and an edge set ${E}$ is a function ${k}:{V} \rightarrow \mathbb{Z}_{>0}$ satisfying ${k(u)}\neq{k(v)}$ for any ${u, v}\in{V}$ such that $\{u, v\}\in{E}$. Let $\mathcal{C}(G)$ be the set of proper colorings of G.

Let $\mathbb{Q}[x_1, \cdots, x_n]$ be the ring of polynomials in $n$ variables with rational coefficients. The symmetric group $\mathfrak{S}_n$ acts on this ring by permuting the variables $x_1,\cdots,x_n$.

\begin{Def}
A polynomial $f$ is \emph{symmetric} if for any permutation $\sigma \in \mathfrak{S}_n$, we have $f(x_1, \cdots, x_n)=f(x_{\sigma(1)}, \cdots, x_{\sigma(n)})$.
\end{Def}

\begin{Def}
For a simple graph $G=({V},{E})$, 
the \emph{chromatic symmetric function} of $G$ is \[X_G(x)=\sum_{k}\prod_{v\in V} x_{k(v)}, \; \text{ for all } k\in\mathcal{C}(G).\]
\end{Def}

\begin{Def}
A polynomial $f$ is \emph{quasisymmetric} if for every composition $\alpha$ of length $k$, the coefficient of $x_1^{\alpha_1}\cdots x_k^{\alpha_k}$ is the same as the coefficient of $x_{i_1}^{\alpha_1}\cdots x_{i_k}^{\alpha_k}$, for any $0<i_1<\cdots<i_k$.
\end{Def}

\begin{Def}
For a simple graph $G=({V},{E})$ with $V\subset\mathbb{Z}_{>0}$, 
the \emph{chromatic quasisymmetric function} of $G$ is \[X_G(x,q)=\sum_{k\in\mathcal{C}(G)} q^{asc(k)}\prod_{v\in V} {x_{k(v)}}_,\] where $asc(k)=\lvert\{\{i,j\}\in E \mid \ i<j\ and \ k(i)<k(j)\}\rvert$.
\end{Def}
Note that the function $X_G(x,1)$ is the chromatic symmetric function.
In general, $X_G(x,q)$ is not a symmetric function. However, it is the symmetric function if $G$ satisfies a certain condition which is explained below.

\begin{Def} [{\cite{shareshian2016chromatic}}]\label{5}
Let $\mathbf{m} \coloneqq (m_1, m_2, \cdots, m_{n-1})$ be a weakly increasing sequence satisfying $i\leq m_i\leq n$ for all $i$. A \emph{natural unit interval order} $P(\mathbf{m})$ is the poset on $[n]$ with the order relation $<_P$ given by $i<_{P(\mathbf{m})}j$ if and only if $m_i<j$.
\end{Def}

\begin{figure}
	\centering
	    \includegraphics[width=\textwidth]{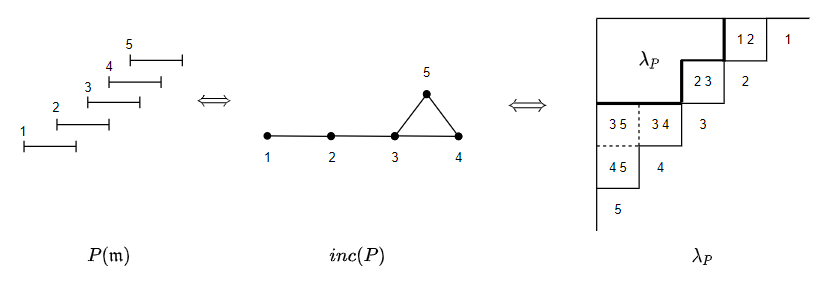}\\
	\caption{$P(\mathbf{m}), inc(P),$ and $\lambda_P$ when $\mathbf{m}=(2,3,5,5)$}
\label{figure1}
\end{figure}

For this paper, it is convenient to work with a partition $\lambda$ defined by $\lambda_i=n-m_i$ for $1\leq i \leq n-1$. Note that there is a one-to-one correspondence between $\mathbf{m}$ and a partition $\lambda$ contained in $\delta_n=(n-1,n-2,\cdots,1)$. By Definition \ref{5}, for a given positive integer $n$ there is a one-to-one correspondence between $\lambda$ and a natural unit interval order $P(\mathbf{m})$, which we denote by $P_\lambda$. (See Figure \ref{figure1})\\

The \emph{incomparability graph} $inc(P)$ of a poset $P$ is a graph defined on the vertex set $P$ so that two elements of $P$ are adjacent if and only if they are incomparable. Now we consider $X_\lambda:=X_{inc(P_\lambda)}$.\\


For a partition $\lambda \subset \delta_n$, define $X_{{\lambda}}(x,q)$ by $X_{P_\lambda}(x,q)$. From this section, we only consider partitions that are contained in $\delta_n$ for a fixed $n$. In this case, we have the following theorem:

\begin{Thm}[{\cite{shareshian2016chromatic}}]\label{thm:symmetric}
For a partition $\lambda \subset \delta_n$, $X_\lambda$ is symmetric.
\end{Thm}

\begin{rem} For a given partition $\lambda\subset \delta_n$, one can notice the corresponding Dyck path. Then area sequence $(a_1,\cdots,a_{n-1})$ of the Dyck path can be computed by $a_i=m_i-i$. For example the area sequence for figure \ref{figure1} is $(1,1,2,1)$.
\end{rem}

\subsection{Conjectures}

Let $\{b_{\lambda}\}$ be a basis of the ring $\Lambda$ of symmetric functions over $\mathbb{Z}$. A function $f$ in $\Lambda$ is \emph{$b$-positive} when $f$ is a linear combination of $\{b_{\lambda}\}$-basis with non-negative coefficients. Note that elementary symmetric functions $\{e_{\lambda}\}$ and complete homogeneous functions $\{h_{\lambda}\}$ are bases of the ring of symmetric functions over $\mb{Z}$, and we call them $e$-basis and $h$-basis respectively.

\begin{Def}
Let $(a_0,\cdots, a_n)$ be a sequence of integers. The sequence is \emph{palindromic} with center of symmetry $\frac{n}{2}$ if $a_j=a_{n-j}$ for $0\leq j \leq n$. The sequence is said to be \emph{unimodal} if
\[a_0\leq a_1 \leq \cdots \leq a_{c-1} \leq a_c \geq a_{c+1}\geq a_{c+2}\geq \cdots \geq a_n\] for some $c$.
\end{Def}

We say the polynomial $P(q):=a_0+a_1 q+\cdots+a_n q^n$ is palindromic and unimodal with center of symmetry $\frac{n}{2}$ if $(a_0, \cdots, a_n)$ has the above properties. We say that $P(q)$ is positive if coefficients $a_i$ are non-negative. For example, a $q$-integer $[n]_q$ is palindromic with center $\frac{n-1}{2}$. Note that the product of palindromic polynomials is still a palindromic polynomial \cite[Corollary 2.3]{sun2015polynomials}.\\

Now we are ready to state conjectures related to the chromatic symmetric functions.

\begin{Def}
A poset $P$ is called \emph{$(r+s)$-free} if it contains no subposet isomorphic to the direct sum of an $r$-element chain and an $s$-element chain.
\end{Def}

For example, the following figure is an example of $(3+1)$-subposet. Therefore, $(3+1)$-free poset means there is no subposet such as figure \ref{figure2}.
\begin{figure}[h]\label{figure2}
	\centering
	    \includegraphics[scale=0.80]{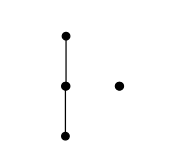}\\
	\caption{Example of $(3+1)$-subposet}
\label{figure2}
\end{figure}

\begin{Con}[Stanley-Stembridge Conjecture {\cite{stanley1993immanants}}]
If $P$ is a (3+1)-free poset, then $X_{inc(P)}(x)$ is $e$-positive.
\end{Con}

Guay-Paquet reduced Stanley's $e$-positivity conjecture for $(3+1)$-free posets to the subclass
of $(3+1)$- and $(2+2)$-free posets  \cite[Theorem 5.1]{guay2013modular}. Such posets are known to be the natural unit interval orders, hence it is enough to consider $X_\lambda(x)$ for $\lambda \subset \delta_n$.\\

Shareshian and Wachs \cite[Corollary 2.8]{shareshian2016chromatic} showed that $X_\lambda(x,q)$ is palindromic, i.e., the coefficients of $X_\lambda(x,q)$, when written as a linear combination of some basis of $\Lambda$, are palindromic with the center of symmetry ${n(n-1)-|\lambda|\over 2}$. They also conjectured the following:

\begin{Con}[{\cite{shareshian2016chromatic}}]\label{nioepositive}$X_\lambda(x,q)$ is $e$-positive and $e$-unimodal. That is,  $X_\lambda(x,q)$ can be written as \[X_\lambda(x,q)=\sum_{\mu}a_{\mu}(q)e_{\mu}(x),\] where the coefficients $a_{\mu}(q)$ are positive and unimodal with center ${n(n-1)-|\lambda|\over 2}$.
\end{Con}
Therefore, Conjecture \ref{nioepositive} states that $X_\lambda(x,q)$ is $e$-positive and $e$-unimodal for any $\lambda \subset \delta_n$.\\

In 1996, Gasharov proved the positivity of Schur function expansion of $X_G(x,q)$ for the incomparability graph of natural unit interval orders when $q=1$  \cite{gasharov1996incomparability}. After 20 years, Shareshain and Wachs proved the same conjecture for general $q$ using $P$-tableaux \cite{shareshian2016chromatic}. They also proved $e$-positivity and $e$-unimodality by using Schur expansion and $P$-tableaux to precisely expand the chromatic quasisymmetric functions with $e$-basis for certain natural unit interval orders.\\

To list some known results, we say that a partition $\lambda \subset \delta_n$ is an \emph{abelian} if it fits inside $\ell \times (n-\ell)$ rectangle for some $\ell$. When $\lambda$ is abelian, there are different formulas for $e$-positivity of $X_\lambda(x,q)$ \cite{cho2019positivity,harada2019cohomology,abrue2020modularlaw}. When $q=1$, the best known result is when the bounce number of $\lambda$ is 3, proved by Cho and Hong in 2019 \cite{cho2020bounce3}.\\

We will relate our work with rook placements and $q$-hit numbers, studied by Abreu and Nigro to describe the coefficients for the $e$-positivity, in Section 4. 

\subsection{Linear relations of chromatic quasisymmetric functions}

Now we are ready to list known linear relations between $X_\lambda$ for different $\lambda$.\\

\begin{Thm}[\cite{lee2018linear} for LLT variant,\cite{guay2013modular} for $q=1$]\label{linearrelation}
Let $\lambda$ be a partition such that $\lambda_i+2\leq\lambda_{i-1}$ for some $i\geq 2$. Let $\mu^0={\lambda}, {\mu}^1, {\mu}^2$ be partitions defined by $\mu^a_j=\lambda_j \text{ if }j\neq i\text{ and }\mu^a_i=\lambda_i+a\text{ for }a=0,1,2$. Then \[X_{{\mu}^0}(x, q)-X_{{\mu}^1}(x, q)=q\big(X_{{\mu}^1}(x, q)-X_{{\mu}^2}(x, q)\big)\] when $\lambda_{n-\lambda_i-1}=\lambda_{n-\lambda_i}$.
\end{Thm}

The following figure \ref{figure3} is an example of $\mu^2, \mu^1,\text{ and } \mu^0$ that satisfies the linear relations from Theorem \ref{linearrelation}.
\begin{figure}[h]\label{figure3}
	\centering
	    \includegraphics[scale=0.7]{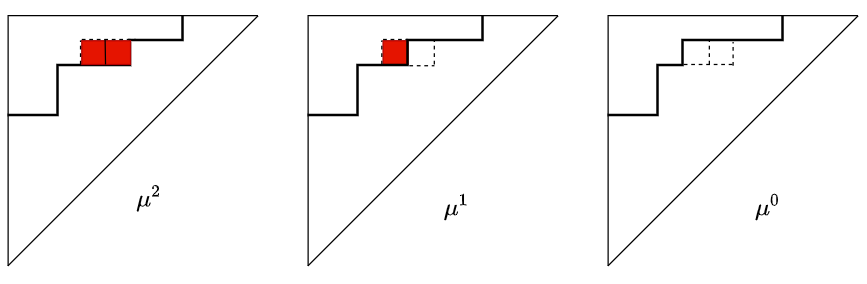}
	\caption{Example of $\mu^2, \mu^1, \text{ and } \mu^0$}
\label{figure3}
\end{figure}

\begin{rem}
Note that the formula in Theorem \ref{linearrelation} is equivalent to
\begin{align}\label{x2x}
X_{{\mu}^0}(x, q)+qX_{{\mu}^2}(x, q)=[2]_qX_{{\mu}^1}(x, q)
\end{align}
When $q=1$, the formula says that $X_{{\mu}^1}(x, q)$ is an average between $X_{{\mu}^0}(x, q)$ and $X_{{\mu}^2}(x, q)$.
\end{rem}

\begin{Thm}[{\cite{huh2020melting}}] \label{linearrelation2}
Let $\lambda$ be a partition such that $\lambda_i+\ell \leq\lambda_{i-1}$ for some $i\geq 2$. Let $\nu^0={\lambda}, {\nu}^1,\cdots, {\nu}^\ell$ be partitions defined by $\nu^a_j=\lambda_j \text{ if }j\neq i\text{ and }\nu^a_i=\lambda_i+a\text{ for }a=0,1,\cdots,\ell$. Assume that $\lambda_{n-\lambda_i}=\lambda_{n-\lambda_i-1}\cdots =\lambda_{n-\lambda_i-\ell+1}$. Then for $0\leq k \leq \ell$, we have

$$[\ell-k]_qX_{\nu^0}(x, q)+q^{\ell-k}[k]_qX_{\nu^\ell}(x, q)=[\ell]_qX_{\nu^k}(x, q)$$

\end{Thm}

The following figure \ref{figure4} is an example of $\nu^0, \nu^\ell, \text{ and } \nu^k$ that satisfies the linear relations from Theorem \ref{linearrelation2}.\\
\begin{figure}[h]
	\centering
	    \includegraphics[scale=0.7]{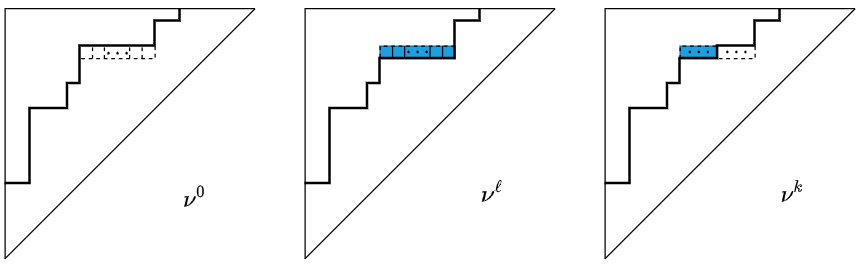}
	\caption{Example of $\nu^0, \nu^\ell, \text{ and } \nu^k$}
\label{figure4}
\end{figure}

Note that there exist column versions of Theorem \ref{linearrelation} and \ref{linearrelation2} by applying $X_\lambda(x,q)=X_{\lambda'}(x,q)$. This identity follows from the definition and their palindromicity.

\section{Explicit formulas for $e$-positivity for abelian cases}

We need the following notations:

\begin{Def}
Let $a_i \in \mathbb{Z}_{\geq0}$ for all $i\in [k]$ and $a=a_1+\cdots+a_k$. Then 
\[{a\choose{a_1,\cdots,a_k}}_q=\frac{[a]_q!}{[a_1]_q! \cdots [a_k]_q!}\text{ where } [a]_q!=[a]_q[a-1]_q\cdots[1]_q\]
is a polynomial in $q$ with non-negative coefficients, which we call \emph{multinomial coefficients} of $q$-integer.
\end{Def}

For this section, we assume that $\lambda$ is abelian, although some of formulas in this section can be applied to non-abelian cases. Let ${\lambda}=(\lambda_1,\cdots,\lambda_\ell)\subset \ell \times (n-s)$ with $\ell\leq s$ and $((n-s)^i)$ be a rectangle $i \times (n-s)$ for $i \leq s$.\\

We introduce explicit formulas of $X_{\lambda}(x,q)$ linearly expanded in terms of $X_{((n-s)^i)}$, $i=0,\cdots,\ell$. Understanding this formula would help to find formulas for $e$-positivity for $X_{\lambda}(x,q)$ as well as the case when $\lambda$ is non-abelian. The following is an example for the explicit formula for $X_{\lambda}(x,q)$ when $\lambda=(\lambda_1,\cdots, \lambda_4),  s=\ell=4$ and $n\geq 8$. For convenience, we denote $X_{\lambda}(x,q)$ by $X_\lambda$.

\begin{align*}
    X_{\lambda}=&\frac{1}{[n-4]_q\cdots[n-7]_q}\times \bigg\{[n-4-\lambda_1]_q\cdots[n-7-\lambda_4]_q X_{((n-4)^0)}\\
    &+\Big(q^{n-7-\lambda_4}[4]_q[\lambda_4]_q[n-6-\lambda_3]_q[n-5-\lambda_2]_q[n-4-\lambda_1]_q\\
    &\quad +q^{n-6-\lambda_3}[3]_q[\lambda_3-\lambda_4]_q[n-5-\lambda_2]_q[n-4-\lambda_1]_q[n-7]_q\\
    &\quad +q^{n-5-\lambda_2}[2]_q[\lambda_2-\lambda_3]_q[n-6-\lambda_4]_q[n-4-\lambda_1]_q[n-7]_q\\
    &\quad +q^{n-4-\lambda_1}[1]_q[\lambda_1-\lambda_2]_q[n-6-\lambda_4]_q[n-5-\lambda_3]_q[n-7]_q\Big)X_{((n-4)^1)}\\
    &+\Big(q^{2n-12-\lambda_3-\lambda_4}(1+q^2)[3]_q[\lambda_4]_q[\lambda_3-1]_q[n-5-\lambda_2]_q[n-4-\lambda_1]_q\\
    &\quad+q^{2n-11-\lambda_2-\lambda_4}[3]_q[2]_q[\lambda_4]_q[\lambda_2-\lambda_3]_q[n-4-\lambda_1]_q[n-7]_q\\
    &\quad+q^{2n-10-\lambda_1-\lambda_4}[3]_q[\lambda_4]_q[\lambda_1-\lambda_2]_q[n-5-\lambda_3]_q[n-7]_q\\
    &\quad+q^{2n-10-\lambda_2-\lambda_3}[3]_q[\lambda_3-\lambda_4]_q[\lambda_2-1]_q[n-4-\lambda_1]_q[n-7]_q\\
    &\quad+q^{2n-9-\lambda_1-\lambda_3}[2]_q[\lambda_3-\lambda_4]_q[\lambda_1-\lambda_2]_q[n-6]_q[n-7]_q\\
    &\quad+q^{2n-8-\lambda_1-\lambda_2}[1]_q[\lambda_2-\lambda_3]_q[\lambda_1-\lambda_4-1]_q[n-6]_q[n-7]_q\Big)X_{((n-4)^2)}\\
    &+\Big(q^{3n-15-\lambda_2-\lambda_3-\lambda_4}[4]_q[\lambda_4]_q[\lambda_3-1]_q[\lambda_2-2]_q[n-4-\lambda_1]_q\\
    &\quad+q^{3n-14-\lambda_1-\lambda_3-\lambda_4}[3]_q[\lambda_4]_q[\lambda_3-1]_q[\lambda_1-\lambda_2]_q[n-7]_q\\
    &\quad+q^{3n-13-\lambda_1-\lambda_2-\lambda_4}[2]_q[\lambda_4]_q[\lambda_2-\lambda_3]_q[\lambda_1-2]_q[n-7]_q\\
    &\quad+q^{3n-12-\lambda_1-\lambda_2-\lambda_3}[1]_q[\lambda_3-\lambda_4]_q[\lambda_2-1]_q[\lambda_1-2]_q[n-7]_q\Big)X_{((n-4)^3)}\\
    &+q^{4n-16-\lambda_1-\lambda_2-\lambda_3-\lambda_4}[\lambda_4]_q[\lambda_3-1]_q[\lambda_2-2]_q[\lambda_1-3]_q X_{((n-4)^4)}\bigg\}
\end{align*}

\begin{figure}[h!] 
	\centering
	    \includegraphics[scale=0.8]{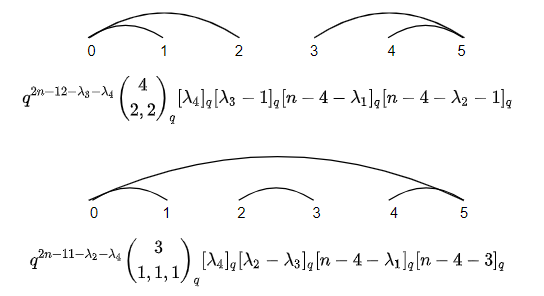}\\
	\centering
		\includegraphics[scale=0.8]{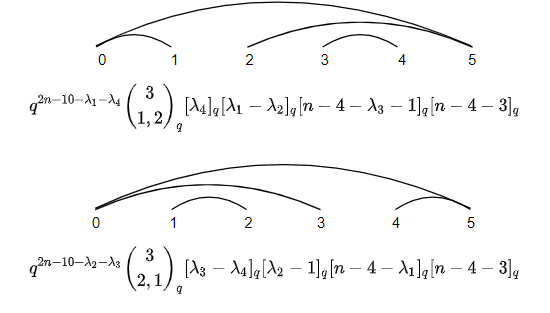}\\
	\centering
		\includegraphics[scale=0.8]{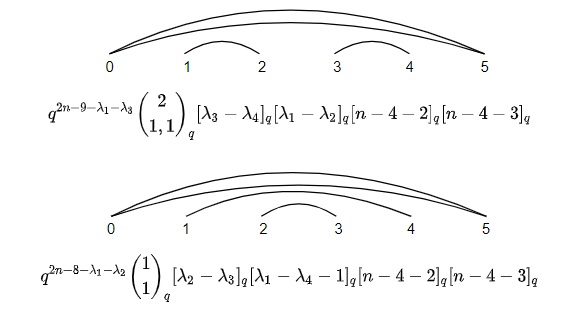}
	\caption{Graphs for coefficients of $X_{((n-4)^2)}$} \label{figure 5}
\end{figure}

Before we describe each term combinatorially, observe that the number of terms in coefficients of $X_{((n-4)^i)}$ is ${4 \choose i}$ and the total number of terms for all coefficients is $\sum_{i=0}^4 {4 \choose i}=2^4$. Graphs in figure \ref{figure 5} are six combinatorial objects that are used to compute coefficients of  $X_{((n-4)^2)}$. (Note that $ {4\choose {2,2}}_q=\frac{[4]_q!}{[2]_q \cdot [2]_q}=[3]_q\cdot \frac{[4]_q}{[2]_q}=[3]_q\cdot (1+q^2)$.)\\

To state the explicit formulas, we first define a graph $G_I$ for a positive integer $r\leq \ell$ and $I\in {[\ell] \choose r}$ as follows. There are $\ell+2$ vertices, labeled by $0,1,\ldots,\ell+1$. For $I=\{ a_1<a_2<\cdots<a_r\}$, join a vertex $a_1$ to a vertex $b_1$ such that the $b_1$ is the largest integer less than $a_1$ contained in $\{0\}\cup I^c$ where $I^c= \{ i \mid 1\leq i \leq \ell, i \notin I\}$.
Second, join a vertex $a_2$ to a vertex $b_2$ such that the $b_2$ is the largest integer less than $a_2$ contained in $\{0\}\cup (I^c\backslash \{b_1\})$. Note that the set always contains $0$.\\

Repeat the procedure for all $i=1,\cdots,r$. Then, join the remaining vertices in $I$ to the vertex $\ell+1$, and then join the vertex $\ell+1$ to the vertex $0$ so that the number of total edges is the same as $\ell$. Note that there could be multiple edges between the vertex $0$ and the vertex $\ell+1$. Also, note that there is no crossing in each graph, i.e., there is no two edges $(p_1,q_1), (p_2,q_2)$ such that $p_1<p_2<q_1<q_2$ by the construction. Note that for a given edge $(p_1,q_1)$, $p_1$ is in $I^c$ or is $0$, and $q_1$ is in $I$ or is $\ell+1$.\\

The following figure \ref{figure6} denotes $G_I$ when $\ell=4, r=2$ and $I=\{1, 4\}$.

\begin{figure}[h!]
	\centering
		\includegraphics[scale=0.85]{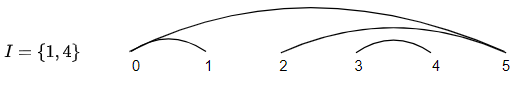}
		\caption{$G_I$ when $\ell=4, r=2$ and $I=\{1, 4\}$}\label{figure6}
\end{figure}

Note that $G_I$ can be described in terms of different combinatorial objects, such as $0-1$ strings, Ferrer diagrams, or link patterns. They appear in the study of maximal parabolic Kazhdan-Lusztig polynomials \cite{Shi2012}. \\

Now we are ready to state the rectangular lemma.

\begin{Thm}[rectangular lemma]\label{rec}
Let ${\lambda}=(\lambda_1,\cdots,\lambda_\ell)\subset \ell \times (n-s)$ with $\ell \leq n-s$. Then $X_\lambda$ is the same as
\[\frac{1}{[n-s]_q[n-s-1]_q\cdots[n-\ell-s+1]_q}\sum_{r=0}^\ell F^{\ell,n-s}_r (\lambda) \cdot X_{((n-s)^r),}\] where
$F^{\ell,n-s}_r(\lambda)$'s are unimodal polynomials in $q$ with non-negative coefficients. Explicitly, we have
\[F^{\ell,m}_r(\lambda)=\sum_{I\in {[\ell] \choose r}}q^{d_{\lambda, I}}c_I \cdot f(\ell,m,{\lambda},I),\]
where the each term appearing in the right-hand side is a non-negative polynomial with the same center.
\end{Thm}    
\begin{rem}
Since all terms in denominators, $[n-s]_q,\cdots,[n-\ell-s+1]_q$, should be greater than $0$, this formula only works when $\ell \leq n-s$. 
\end{rem}

Before we prove rectangular lemma, we first define $f(s,m, {\lambda},I), c_I,$ and ${d_{\lambda, I}}$ in order.
First, we define $f(\ell, m, {\lambda},I)$ associated to each $G_I$ as follows:\\

Let $\lambda_i=0$ if $i>s$ for convenience. Define $\lambda_0$ by $m$, the second input for the function $f$. For the proof of Theorem \ref{rec}, $m$ is equal to $n-s$. Define
    \[f(\ell,m, {\lambda},I):=\prod [\lambda_{\ell+1-a_i}-\lambda_{\ell+1-b_i}-l_i+1]_q,\] 
    where the product runs over all edges $(b_i,a_i)$ of $G_I$ and $l_i$ is the number of edges $(b_j,a_j)$ such that $b_i\leq b_j < a_j \leq  a_i$. Note that $j$ can be equal to $i$ hence $l_i\geq 1$. We call the number \emph{the length of the edge} $(b_i,a_i)$, denoted by $|(b_i,a_i)|$. In other word, the length of an edge is defined by the number of edges under the given edge when we draw edges similar to the above figures. If there are $x$ multiple edges $(0, \ell+1)$, then define the lengths of those $x$ edges by $\ell-x+1,\cdots,\ell$. \\
 
Although $f(\ell,m,\lambda,I)$ is well-defined for any integer $m$, it is not always positive. For Theorem \ref{rec}, we have $m=n-s \geq \lambda_1$ so that $(m,\lambda_1,\cdots,\lambda_\ell)$ is a partition. Then the following holds:

\begin{Lem} \label{fpositive}
If $f(\ell,n-s,\lambda,I)$ is not zero, then each term in $f(\ell,n-s,\lambda,I)$ is a positive polynomial. 
\end{Lem}
\begin{proof}
We will show that if $f$ is nonzero, each term $[\lambda_{\ell+1-a_i}-\lambda_{\ell+1-b_i}-l_i+1]_q$ is a positive polynomial. 
Assume otherwise. Consider a pair $(b_i,a_i)$ such that $\lambda_{\ell+1-a_i}-\lambda_{\ell+1-b_i}-l_i+1$ is negative and $l_i$ is minimal. If $l_i$ is 1 then the term is obviously positive. If $l_i>1$, then there exist a positive integer $x$ and integers $b_i+1= p_1<p_2<\cdots<p_x<p_{x+1}=a_i$ where $(p_1,p_2-1),(p_2,p_3-1),\ldots,(p_x,p_{x+1}-1)$ are edges of $G_I$. By the construction, $ \lambda_{\ell+2-p_{y+1}}-\lambda_{\ell-p_{y}+1}-l_i+1$ is positive, i.e.,  $ \lambda_{\ell+2-p_{y+1}}-\lambda_{\ell+1-p_{y}}-l_i\geq 0$. Let $l^{(j)}$ be the length of $(p_{j},p_{j+1}-1)$. Then $l_i=1+\sum_{j=1}^x l_i^{(j)}$ and 
\begin{align*}
\lambda_{\ell+1-a_i}-\lambda_{\ell+1-b_i}-l_i+1 & \geq \lambda_{\ell+2-p_{x+1}}-\lambda_{\ell+1-p_{1}}-l_i+1\\
&=\sum_{j=1}^x (\lambda_{\ell+2-p_{j+1}}-\lambda_{\ell+1-p_j}-l^{(j)})\\
&\geq 0,
\end{align*}
which makes a contradiction.
\end{proof}

The following figure is a graph for $f(\ell,n-s,{\lambda}, I)$ when ${\lambda}=(\lambda_1,\cdots,\lambda_4), \ell=s=4, \;r=2, \; I=\{1,4\}$ and $(b_1,a_1)=(0,1), (b_2,a_2)=(3, 4), (b_3,a_3)=(2, 5), (b_4,a_4)=(0,5)$.

\begin{figure}[h!] 
	\centering
		\includegraphics[scale=0.85]{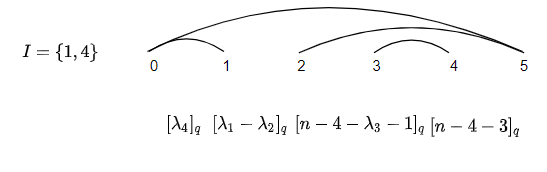}
		\caption{$f(4,n-4,\lambda, I)$ when $I=\{1, 4\}$}
\label{114} 
\end{figure}

Secondly, we define $c_I$ which only depends on $I\in {[\ell] \choose r}$. Since the set $I$ determines the graph $G$, we define $c_I$ for a given graph $G$ and we denote by $c_G$. First of all, if the vertex $0$ and $\ell+1$ are not connected, add an edge between them. If there are multiple edges between the vertex $0$ and $\ell+1$, delete all edges except one edge. Now we define the $c_G$ inductively. At each step, we delete the edge of the longest length. Then there exist integers $0\leq p_1<p_2<\cdots<p_y<p_{y+1}\leq \ell +2$ where the vertex $p_1$ is the non-isolated vertex of minimal label (which should be either $0$ or $1$), $p_{y+1}-1$ is the non-isolated vertex of maximal label, and $(p_1,p_2-1),(p_2,p_3-1),\ldots,(p_y,p_{y+1}-1)$ are edges of $G$. For example, when $G$ is the graph in Figure \ref{114}, we have $p_1=0$, $p_2=2$, $p_3=6$.

For $i$-th edge $(p_{i},p_{i+1}-1)$ of such edges, we can define a induced subgraph $G_i$ of $G$ with vertices $p_{i},p_{i}+1,\ldots,p_{i+1}-1$. Let $n_i$ be the number of edges of $G_i$. Then $c_G$ is inductively defined by
\[{\sum_{i=1}^y n_i \choose n_1,\cdots,n_s}_q\;\prod_{i=1}^y c_{G_i}.\]
Note that to compute $c_{G_i}$ for each $i$, one needs to relabel the vertices of $G_i$ so that the left most vertex $p_i$ is relabeled by $0$, the vertex $p_i+1$ is relabeled by $1$, etc. Also for a graph $G'$ with one edge, define $c_{G'}$ by $1$. For the graph in Figure \ref{114}, we have $c_G={3 \choose 1,2}_q c_{G_1}c_{G_2}={3 \choose 1,2}_q=[3]_q$. The above formula shows that $c_G$ is in $\mathbb{N}[q]$, and by inductively using the definition of $c_H$ for all $H$ appearing during the induction, we have the following:
\begin{Lem} \label{cproduct}
We have
$$ c_G= { [\ell]_q! \over \prod_{(b_i,a_i)} [|(b_i,a_i)|]_q }$$
where $(b_i,a_i)$ runs over all $\ell$ edges of $G$. 
\end{Lem}
The proof of Lemma \ref{cproduct} directly follows from the inductive definition of $c_G$. Note that it is easier to compute $c_G$ but it is not clear that the right-hand side is indeed in $\mathbb{N}[q]$. \\
\begin{rem}
The term $[\ell]_q! \over \prod_{(b_i,a_i)} [|(b_i,a_i)|]_q $ appears in a $q$-analog of Knuth's tree hook-length formula studied by Bjorner and Wachs \cite{BW1989}, hence the definition of $c_G$ in this paper shows another proof that $c_G$ is in $\mathbb{N}[q]$.
\end{rem}

At last, for $I=\{a_1, \cdots, a_r\}$, define $d_{\lambda, I}$ as follows:
\[d_{\lambda, I}:= r(n-s)-\sum_{i=1}^{r}(\ell-r+i-a_i)-\sum_{i=1}^{r}\lambda_{\ell-a_i+1.}\]
    
We first show the following:

\begin{Thm} All terms $q^{d_{\lambda, I}}c_I \cdot f(\ell,n-s,{\lambda}, I)$ in $F_r^{\ell,n-s}(\lambda)$ have the same center when $r$ is fixed.

\begin{proof}
Let $\lambda = (\lambda_1, \cdots, \lambda_{\ell})$ and let $r$ is fixed. Then the degree of $q^{d_{\lambda, I}}$ is \[r(n-s)-\sum_{i=1}^{r}(\ell-r+i-a_i)-\sum_{i=1}^{r}(\lambda_{\ell-a_i+1}),\]

and the center of  $c_I \text{ and } f(\ell,n-s,{\lambda},I)$ are
\[\frac{1}{2}(\frac{\ell(\ell+1)}{2}-\sum_{i=1}^{\ell}l_i) \text{ and } \frac{1}{2}\sum_{i=1}^{\ell}(\lambda_{\ell-a_i+1}-\lambda_{\ell-b_i+1}-l_i), \text{respectively.}\]

In addition, note that \[\frac{1}{2}\sum_{i=1}^{\ell}(\lambda_{\ell-a_i+1}-\lambda_{\ell-b_i+1}) - \sum_{i=1}^{r}(\lambda_{\ell-a_i+1}) = -\frac{1}{2}\{\lambda_1 + \cdots + \lambda_{\ell}\}+\frac{1}{2}\{ (n-s)(\ell-r)\}\]

Thus, the center of $q^{d_{\lambda, I}}c_I \cdot f(\ell,n-s,{\lambda}, I)$ is as follows :

\begin{align*}
    &\text{\scalebox{0.83}{$\{r(n-s)-\sum_{i=1}^{r}(\ell-r+i-a_i)-\sum_{i=1}^{r}(\lambda_{\ell-a_i+1})\}+\frac{1}{2}\{\sum_{i=1}^{\ell}(\lambda_{\ell-a_i+1}-\lambda_{\ell-b_i+1}-l_i)+(\frac{\ell(\ell+1)}{2}-\sum_{i=1}^{\ell}l_i)\}$}}\\
    &=\text{\scalebox{0.9}{$r(n-s)-\sum_{i=1}^{r}(\ell-r+i-a_i)+\frac{1}{2}\{\{\frac{\ell(\ell+1)}{2}-2\sum_{i=1}^{\ell}l_i\}-(\lambda_1+\cdots+\lambda_{\ell})+(n-s)(\ell-r)\}$}}\\
    &=\text{\scalebox{0.9}{$r(n-s)-\sum_{i=1}^{r}(\ell-r+i)+\sum_{i=1}^{r}a_i-\sum_{i=1}^{\ell}l_i+\frac{1}{2}\{\frac{\ell(\ell+1)}{2}-(\lambda_1+\cdots+\lambda_{\ell})+(n-s)(\ell-r)\}$}}\\
    &=\text{\scalebox{0.83}{$r(n-s)-r(\ell-r)-\frac{r(r+1)}{2}-\frac{(\ell-r)(\ell-r+1)}{2}+\frac{1}{2}\{\frac{\ell(\ell+1)}{2}-(\lambda_1+\cdots+\lambda_{\ell})+(n-s)(\ell-r)\}$}}\\
 &=\text{\scalebox{0.83}{${(n-s)(\ell+r) \over 2} - {\ell(\ell+1)\over 4}-{|\lambda| \over 2}$}}
\end{align*}
From the third line to the fourth line, we apply the identity
$$\sum_{i=1}^{r}a_i-\sum_{i=1}^{\ell}l_i=-\frac{(\ell-r)(\ell-r+1)}{2},$$
which can be proved by an induction on $\sum_{i=1}^r{a_i}$ with the initial case $I=\{1,2,\cdots,r\}$. Therefore, since there is no terms with $\{a_i, l_j\}$, the center of each term in $F_r^{\ell,n-s}(\lambda)$ is the same when $r$ is fixed. One can also check that the center of $F_r^{\ell,n-s}(\lambda)$ does match with the equation in Theorem \ref{rec}, namely, the center of $F_r^{\ell,n-s}(\lambda) X_{(n-s)^r}$ is the same as the center of $X_\lambda$ plus the center of $[n-s]_q \cdots [n-s-\ell+1]_q$.
\end{proof}
\end{Thm}

Before we prove Theorem \ref{rec}, the following lemma is useful:
\begin{Lem} 

\begin{enumerate}[(a)]\label{powerlemma}
\item For given $\lambda$ and $I$, if $\lambda_{\ell-j}=\lambda_{\ell-j+1}$, $j \notin I$ and $j+1\in I$,
then we have 
$$f(\ell,m,\lambda,I)=0.$$
\item If $f(\ell,m,\lambda,I)$ is nonzero and $\lambda_{\ell-j-x+1}=\cdots=\lambda_{\ell-j+1}$ for some $1\leq x \leq \ell-j+1$, then there exists the unique $0\leq y\leq x$ such that
\begin{align}\label{iint} 
I \cap \{j,\cdots,j+x\}= \{j,\cdots,j+y\}.
\end{align}
Here, $y=0$ means that $\{j, j+1, \cdots,j+y\}$ is the empty set.
\end{enumerate}

\end{Lem}
\begin{rem}
For Lemma \ref{powerlemma}, $m$ can be any integer, namely, $m\geq \lambda_1$ is not necessary.
\end{rem}
The first lemma directly follows since $(j,j+1)$ is an edge in $G_I$. The lemma suggests that to get a nonzero term, when choosing $I$, one has to choose $j$ before $j+1$ when $\lambda_{\ell-j}=\lambda_{\ell-j+1}$. The second lemma also follows since if Equation (\ref{iint}) does not hold, there exists a number $a \in [j,j+x-1]$ such that $a\notin I$ and $a+1\in I$. This makes a contradiction.\\

Now we prove the rectangular lemma.
\begin{proof}[Proof of rectangular lemma]
We need three steps to prove the theorem \ref{rec}: 
\begin{enumerate}[Step I]
\item Show that the theorem holds when $\lambda=((n-s)^{s'})$ for all $0\leq s'\leq \ell$.
\item Show that $F_r^{\ell,n-s}(\lambda)$ satisfies the row and column linear relations in Theorem \ref{linearrelation}.
\item Prove that $F_r^{\ell,n-s}((n-s)^{s'})$ from the first step and the linear relations from second step determine $F_r^{\ell,n-s}(\lambda)$ for any $\lambda\subset \ell \times (n-s)$.
\end{enumerate}

\textbf{Step I}. \textit{The Theorem \ref{rec} holds when $\lambda'=((n-s)^{s'})$ for all $0\leq s'\leq \ell$.}

\begin{proof}
First of all, we need to show $F_r^{\ell,n-s}((n-s)^{s'})=0$ for $s'\neq r$. We first show that $f(\ell,n-s,(n-s)^{s'},I)=0$ for all $I$ when $s^{\prime}\neq r$, and $[n-s]_q\cdots[n-s-\ell+1]_q$ if $s'=r$. \\

Assume that $f(\ell,n-s,(n-s)^{s'},I)$ is nonzero. Then we have $\ell \in I$ because otherwise $(\ell,\ell+1)$ is an edge. Then by Lemma \ref{powerlemma}, $I$ must be $\{\ell-r+1,\ell-r+2,\cdots, \ell\}$.
Therefore, $s'=r$ and $F_{(n-s)^{r}}(r)=[n-s]_q[n-s-1]_q\cdots[n-s-\ell+1]_q$, since $d_{(n-s)^r,I}=c_I=1$ for $I=\{\ell-r+1,\ell-r+2,\cdots, \ell\}$. It implies that 
\begin{align*}
    &\frac{1}{[n-s]_q[n-s-1]_q\cdots[n-\ell-s+1]_q}\sum_{s'=0}^\ell F^{\ell,n-s}_{s'} (\lambda) \cdot X_{((n-s)^{s'}),} \\
  &=\frac{1}{[n-s]_q\cdots[n-s-\ell+1]_q}\times[n-s]_q[n-s-1]_q\cdots[n-s-\ell+1]_q X_{((n-s)^r)}\\
  &  =X_{((n-s)^r)}
\end{align*}
\end{proof}
\textbf{Step II}. \textit{$F_\lambda(r)$ satisfies the row and column linear relations in Theorem \ref{linearrelation}.}
\begin{proof}
We have to check that the stated formula satisfies the linear relations from Theorem \ref{linearrelation}. We divide the linear relations into the row relation and the column relation. We show the row relation first.\\

Since $s,\ell,n$ does not vary in the proof, we simply write $f(\ell,n-s,\lambda,I)$ by $f(\lambda,I)$, and $F_r^{\ell,n-s}(\lambda)$ by $F_r(\lambda)$. Let $\lambda^0={\mu}, {\lambda}^1, {\lambda}^2$ be partitions defined by $\lambda^a_j=\mu_j \text{ if }j\neq i\text{ and }\lambda^a_i=\mu_i+a\text{ for }a=0,1,2$. We have to show $F_{r}(\lambda^0)+q\cdot F_{r}(\lambda^2)=[2]_q F_{r}(\lambda^1)$. Since $c_I$ only depends on $I$, thus when $r$ and $I$ are fixed, it is enough to show 
\begin{align}\label{3} q^{d_{0}}f(\lambda^0,I)+q\cdot q^{d_{2}}f(\lambda^2,I)=[2]_q q^{d_{1}}f(\lambda^1,I)\end{align}
where $d_i=d_{\lambda^i, I}$. In addition, since $f(\lambda^j, I)$ for $ j=0,1,2$ have the same terms except the term containing $\mu_i$, it is enough to consider the terms with $\mu_i$.

If $\ell-i+1 \in I$, then, by the way $d_{\lambda, I}$ defined, $d_0=d_1+1=d_2+2$. Then we have $q^{d_0}[\mu_i-K]_q+q\cdot q^{d_2}[\mu_i+2-K]_q=q^2\cdot q^{d_2}[\mu_i-K]_q+q\cdot q^{d_2}[\mu_i+2-K]_q=q\cdot q^{d_2} [2]_q[\mu_i+1-K]_q=q^{d_1} [2]_q[\mu_i+1-K]_q$ for any integer $K$. Therefore, (\ref{3}) holds.\\

If $\ell-i+1 \notin I$, then $d_0=d_1=d_2$. Thus we have $ q^{d_0}[M-\mu_i]_q+q\cdot q^{d_2}[M-\mu_i-2]_q=q^{d_1}[2]_q[M-\mu_i-1]_q$ for any integer $M$. Therefore, the row relation from Theorem \ref{linearrelation} holds for $F_r(\lambda)$.\\

Next, we need to show this formula satisfies the column relation. For a partition $\mu$ with $i$ such that $\mu_{i+1}=\mu_i\leq \mu_{i-1}-1$, let $\lambda^{(0)}=\mu, \lambda^{(1)}, \lambda^{(2)}$ be partitions defined by $\lambda_j^{(1)}=\mu_j$ if $j\neq i$, $\lambda_i^{(1)}=\mu_i+1$,\; and $\lambda_j^{(2)}=\mu_j$ if $j\neq i,i+1$, $\lambda_i^{(2)}=\mu_i+1=\lambda_{i+1}^{(2)}$. We want to show $F_{\lambda^{(0)}}(r)+q\cdot F_{\lambda^{(2)}}(r)=[2]_q F_{\lambda^{(1)}}(r)$.\\
 
For $X_{\lambda^{(1)}}$, it can be divided into 2 parts by $I$, $A$ and $B$, where $A$ is the collection of all $I$ such that $\ell-i+1\in I$, $\ell-i\notin I$ and $B$ is all other $I$'s. Note that if $I\in A$, then $f(\lambda^{(0)},I)=f(\lambda^{(2)},I)=0$. Thus, it is enough to show 
\begin{align}
   \text{\scalebox{0.98}{$[2]_q\sum_{\substack{I\in {[\ell] \choose r}\\ I\in A}}q^{d_{1, I}}c_{I} f_1(I) =\sum_{\substack{H\in {[\ell] \choose r}\\ H \in B}}c_{H}\Big(q^{d_{0, H}} f_0(H)+q\cdot q^{d_{2, H}} f_2(H)-[2]_q q^{d_{1, H}} f_1(H)\Big)_,$}} \label{totalsum}
\end{align}
where $f_i(I)=f(\lambda^{(i)},I)$ and $d_{i, I}=d_{\lambda^{(i)}, I}$ for $i=0,1,2$. Brief idea for proving the above equation is as follows: For a given $I'$ in $A$, consider the subset $B'\subset B$ consisting of a graph $H'$ such that $G_{H'}$ can be obtained from $G_{I'}$ by deleting two edges including $(\ell-i,\ell-i+1)$ and then add two edges so that the resulting graph is of the form $G_{H'}$. Note that there should be no crossing in $G_{H'}$.\\

To describe all possible $H'$, fix $I'\in A$. Then there exists the unique edge $\delta$ (denote by $(\delta_b, \delta_a)$) in $G_{I'}$ such that the edge is the shortest edge with $\delta_b<\ell-i<\ell-i+1<\delta_a$. 
Then there exist edges $(\beta_{1b},\beta_{1a}),(\beta_{2b},\beta_{2a}),\cdots,(\beta_{m_2b},\beta_{m_2a})$ of $G_{I'}$ such that $\beta_{1b}=\ell-i+2, \beta_{(j'+1)b}=\beta_{j'a}+1$ for $j'=1,2,\cdots m_2-1$, and $\beta_{m_2a}=\delta_a-1$. Denote $(\beta_{j'b},\beta_{j'a})$ by $\beta_{j'}$. In other words, edges $\beta_{j'}$ are consecutive edges between the vertices $\ell-i+2$ and $\delta_a$. Similarly, define edges $\alpha_j=(\alpha_{jb},\alpha_{ja})$ by $\alpha_{1a}=\ell-i-1, \alpha_{(j+1)a}=\alpha_{jb}-1$ for $j=1, 2, \cdots, m_1-1$,  and $\alpha_{m_1b}=\delta_b+1$, therefore, edges $\alpha_j$ are consecutive edges between vertices $\delta_b$ and $\ell-i-1$. See Figure \ref{column1}.

\begin{figure}[h!]
	\centering
		\includegraphics[width=\textwidth]{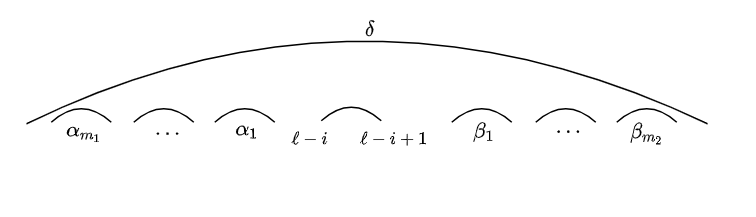}
		\caption{The graph $G_I$ under the edge $\delta$}\label{column1}
\end{figure}

Let $G_{\alpha_j}$ be the graph associated with the graph $G_I$ such that the edge $(\alpha_{jb}, \alpha_{ja})$ and $(\ell-i, \ell-i+1)$ are switched into $(\alpha_{ja}, \ell-i)$ and $(\alpha_{jb}, \ell-i+1)$, respectively. Define the graph $G_{\beta_{j^{\prime}}}$ similarly. Let $G_{\delta}$ be the graph obtained from $G_I$ by replacing edges $\delta$ and $(\ell-i,\ell-i+1)$ by $(\delta_b,\ell-i)$ and $(\ell-i+1,\delta_a)$. Then it is clear that the subset $B' \subset B$ consists of
\begin{enumerate}[(a)]
\item $I'/\{\alpha_{ja}\} \cup \{ \ell-i \}$ for $j=1,2,\ldots, m_1$
\item $I'/\{\ell-i+1\}\cup\{ \beta_{j^{\prime}b}\}$ for $j^{\prime}=1,2,\ldots,m_2$
\item $I'/\{\ell-i+1\}\cup \{ \ell-i\}$.
\end{enumerate} 

Denote the set in (a) by $I_{\alpha_j}$, the set in (b) by $I_{\beta_{j^{\prime}}}$, and the set (c) by $I_\delta$. For simplicity, We often denote $G_{I_{\alpha_j}}$ by $G_{\alpha_j}$, and we do similarly for $\beta_{j^{\prime}}$ and $\delta$. Figure \ref{column2} is an example of a graph $G_{\alpha_j}$ associated with graph $G_I$.\\

Let $[G^k_{I'}]_q$ be $f(\lambda^{(k)},I')$ and let $[G^k_{e}]_q$ be $f(\lambda^{(k)},I_e)$ for an edge $e \in \{ \delta,\alpha_j,\beta_{j^{\prime}} \mid j=1,2,\cdots,m_1, j^{\prime}=1,2,\cdots,m_2\}$ when $k=0,1,2$. Let $d_{k,J}$ be $d_{\lambda^{(k)}, J}$ when $k=0,1,2$ and $J\in \{I,I_{\alpha_j},I_{\beta_{j^{\prime}}},I_\delta \mid j=1,\cdots, m_1$, $j^{\prime}=1,\cdots,m_2\}$. \\
\begin{figure}[h!]
	\centering
		\includegraphics[width=\textwidth]{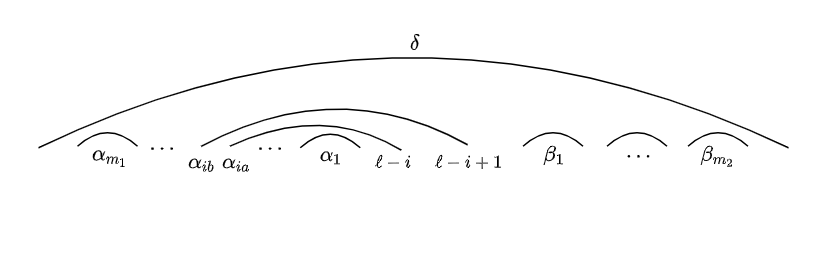}
		\caption{The graph $G_{\alpha_j}$ associated with the graph $G_I$} \label{column2}
\end{figure}

Now we are ready to prove Equation (\ref{totalsum}).
\begin{Pro}\label{Pro1}
For given $j\leq m_1$ and $j^{\prime}\leq m_2$, we have the following:
\begin{enumerate}[(a)]
\item $q^{d_{1, I}+1+\sum_{l=1}^{j-1}\lvert \alpha_l \rvert}[G^1_{I}]_q=q^{d_{0, {I}_{\alpha_j}}}[G^0_{\alpha_j}]_q+q\cdot q^{d_{2, {I}_{\alpha_j}}}[G^2_{\alpha_j}]_q-[2]_q q^{d_{1, {I}_{\alpha_j}}}[G^1_{\alpha_j}]_q.$
\item $q^{d_{1, I}+1+\sum_{l=1}^{j^{\prime}-1}\lvert \beta_l \rvert}[G^1_{I}]_q=q^{d_{0, {I}_{\beta_{j^{\prime}}}}}[G^0_{\beta_{j^{\prime}}}]_q+q\cdot q^{d_{2, {I}_{\beta_{j^{\prime}}}}}[G^2_{\beta_{j^{\prime}}}]_q-[2]_q q^{d_{1, {I}_{\beta_{j^{\prime}}}}}[G^1_{\beta_{j^{\prime}}}]_q.$
\item $q^{d_{1,I}}[G^1_{I}]_q=q^{d_{0, {I}_{\delta}}}[G^0_{\delta}]_q+q\cdot q^{d_{2, I_{\delta}}}[G^2_{\delta}]_q-[2]_q q^{d_{1, I_{\delta}}}[G^1_{\delta}]_q.$
\end{enumerate}
\end{Pro}

\begin{proof} We will only prove $(a)$ and $(c)$, since the proof of $(b)$ is similar to $(a)$. \\
Note that $d_{0, I_{\alpha_j}}=1+d_{1, {I}_{\alpha_j}}=2+d_{2, I_{\alpha_j}}$. Also, we have 
\begin{align*}
d_{1, I} &={d_{1, I_{\alpha_j}}}+{\mu_i-\mu_{\ell-\alpha_{ja}+1}-(\ell-i)+\alpha_{ja}}\\
&={d_{1, {I}_{\alpha_j}}}+{\mu_i-\mu_{\ell-\alpha_{ja}+1}-2\sum_{l=1}^{l=j-1}\lvert \alpha_l \rvert-1}
\end{align*}
Also note that $[G^1_{I}]_q$ and $[G^k_{\alpha_j}]_q$ have exactly the same terms except terms with $\mu_i, \mu_{\ell-\alpha_{ja}+1}$, and $\mu_{\ell-\alpha_{jb}+1}$ for $k=0,1,2$. Therefore, it is enough to check terms with $\mu_i, \mu_{\ell-\alpha_{ja}+1}$, and $\mu_{\ell-\alpha_{jb}+1}$. We have
\begin{align*}
RHS \text{ of }(a)&=q^{d_{0, {I}_{\alpha_j}}}[G^0_{\alpha_j}]_q+q\cdot q^{d_{2, I_{\alpha_j}}}[G^2_{\alpha_j}]_q-[2]_q q^{d_{1, {I}_{\alpha_j}}}[G^1_{\alpha_j}]_q\\
 &=q^{d_{1, {I}_{\alpha_j}}}\Big(q\cdot\big[\mu_i-\mu_{\ell-\alpha_{ja}+1}-\sum_{l=1}^{j-1}\lvert \alpha_l \rvert\big]_q\cdot\big[\mu_i-\mu_{\ell-\alpha_{jb}+1}-\sum_{l=1}^{j}\lvert \alpha_l \rvert\big]_q\\
    &+\big[\mu_i+1-\mu_{\ell-\alpha_{ja}+1}-\sum_{l=1}^{j-1}\lvert \alpha_l \rvert\big]_q\cdot\big[\mu_i+1-\mu_{\ell-\alpha_{jb}+1}-\sum_{l=1}^{j}\lvert \alpha_l \rvert\big]_q\\
    &-[2]_q\big[\mu_i-\mu_{\ell-\alpha_{ja}+1}-\sum_{l=1}^{j-1}\lvert \alpha_l \rvert\big]_q\cdot\big[\mu_i+1-\mu_{\ell-\alpha_{jb}+1}-\sum_{l=1}^{j}\lvert \alpha_l \rvert\big]_q\Big)\\
    &=q^{d_{1, {I}_{\alpha_j}}}\Big(\big[\mu_i+1-\mu_{\ell-\alpha_{jb}+1}-\sum_{l=1}^{j}\lvert \alpha_l \rvert\big]_q-\big[\mu_i      -\mu_{\ell-\alpha_{ja}+1}-\sum_{l=1}^{j-1}\lvert \alpha_l \rvert\big]_q\Big)\\
    &=q^{d_{1, {I}_{\alpha_j}}}\cdot q^{\mu_i-\mu_{\ell-\alpha_{ja}+1}-\sum_{l=1}^{j-1}\lvert \alpha_l \rvert}\big[\mu_{\ell-\alpha_{ja}+1}-\mu_{\ell-\alpha_{jb}+1}-(\lvert\alpha_j\rvert-1)\big]_q\\
&=q^{d_{1,I}+1+\sum_{l=1}^{j-1}\lvert \alpha_l \rvert}\big[\mu_{\ell-\alpha_{ja}+1}-\mu_{\ell-\alpha_{jb}+1}-(\lvert\alpha_j\rvert-1)\big]_q=LHS \text{ of } (a).
\end{align*}

For proving Part $(c)$, first note that 
$$d_{1,I}=d_{0,I_{\delta}}=d_{1,I_{\delta}}=d_{2,I_{\delta}}+1.$$
Then after canceling all common terms, we have
\begin{align*}
&q^{d_{0, {I}_{\delta}}} [\mu_i-K]_q[M-\mu_i]_q+ q^{1+d_{2, I_{\delta}}} [1+\mu_i-K]_q[M-\mu_i-1]_q \\
&-[2]_q q^{d_{1, I_{\delta}}} [\mu_i-K]_q[M-\mu_i-1]_q.\\
&=q^{d^{1,I}} [M-K-1]_q
\end{align*}
where $K=\mu_{\delta_b}-\sum_{j=1}^{m_1}\lvert \alpha_j \rvert$ and $M=\sum_{j^{\prime}=1}^{m_2}\lvert \beta_{j^{\prime}} \rvert +\mu_{\delta_a}$. Note that, 
$$M-K-1= \mu_{\delta_a}- \mu_{\delta_b} - \lvert \delta\rvert +1$$
hence we proved (c).
\end{proof}

Note that we have not used the term $c_I$. The following proposition related to $c_I$ is enough to show that $F_\lambda(r)$ satisfies the column linear relation.

\begin{Pro}\label{Pro2} For a given $I$ satisfying $\ell-i \notin I$ and $\ell-i+1 \in I$, we have
\[[2]_qc_{{G}_I}=\sum_{j=1}^{m_1}q^{1+\sum_{l=1}^{j-1}\lvert \alpha_l \rvert}c_{G_{I_{\alpha_j}}}+\sum_{j^{\prime}=1}^{m_2}q^{1+\sum_{l=1}^{j^{\prime}-1}\lvert \beta_l \rvert}c_{G_{I_{\beta_{j^{\prime}}}}}+c_{G_{I_{\delta}}}.\]\\
\end{Pro}

Note that Propositions \ref{Pro1} and \ref{Pro2} prove the column linear relation of $F_r(\lambda)$. By Lemma \ref{cproduct}, one can compute the following ratios:
\begin{align*}
{c_{G_{I_{\alpha_j}}} \over c_{G_I}}= { [\lvert \alpha_j\rvert]_q \over [1+ \sum_{l=1}^{j-1} \lvert \alpha_l \rvert ]_q \cdot [1+ \sum_{l=1}^j \lvert \alpha_l \rvert ]_q}.\\
{c_{G_{I_{\beta_{j^{\prime}}}}} \over c_{G_I}}= { [\lvert \beta_{j^{\prime}}\rvert]_q \over [1+ \sum_{l=1}^{{j^{\prime}}-1} \lvert \beta_l \rvert ]_q \cdot [1+ \sum_{l=1}^{j^{\prime}} \lvert \beta_l \rvert ]_q}.\\
{c_{G_{I_{\delta}}} \over c_{G_I}}= { [2+\sum_{l=1}^{m_1} \lvert \alpha_l \rvert+\sum_{l=1}^{m_2} \lvert \beta_l \rvert]_q \over [1+ \sum_{l=1}^{m_1} \lvert \alpha_l \rvert ]_q \cdot [1+ \sum_{l=1}^{m_2} \lvert \beta_l \rvert ]_q}.
\end{align*}

By the following two lemmas, the proposition \ref{Pro2} is proved. 

\begin{Lem}\label{Lem1}
Let $a_i\in \mathbb{Z}_{\geq0}$ for $i=1,\cdots,n$, and $a=\sum_{i=1}^n a_i$. Then we have
\begin{align*}
    \sum_{i=1}^n q^{a_1+\cdots+a_{i-1}+1} \frac{[a_i]_q}{[a_1+\cdots+a_{i-1}+1]_q[a_1+\cdots+a_{i}+1]_q}=1-\frac{1}{[a+1]_q}_.
\end{align*}
\begin{proof}[Proof of Lemma \ref{Lem1}]
    For all $i=1,\cdots,n$, \[q^{a_1+\cdots+a_{i-1}+1}[a_i]_q=[a_1+\cdots+a_{i}+1]_q-[a_1+\cdots+a_{i-1}+1]_q.\]
    Therefore,
    \begin{align*}
        &\sum_{i=1}^n q^{a_1+\cdots+a_{i-1}+1} \frac{[a_i]_q}{[a_1+\cdots+a_{i-1}+1]_q[a_1+\cdots+a_{i}+1]_q}\\
        &=\sum_{i=1}^n\frac{[a_1+\cdots+a_{i}+1]_q-[a_1+\cdots+a_{i-1}+1]_q}{[a_1+\cdots+a_{i-1}+1]_q[a_1+\cdots+a_{i}+1]_q}\\
        &=\sum_{i=1}^n \frac{1}{[a_1+\cdots+a_{i-1}+1]_q}-\frac{1}{[a_1+\cdots+a_{i}+1]_q}=1-\frac{1}{[a+1]_q}_.
    \end{align*}
\end{proof}
\end{Lem}

\begin{Lem}\label{Lem2}
Let $a, b\in \mathbb{Z}_{\geq0}$. Then we have
    \begin{align*}
        \frac{[a+b+2]_q}{[a+1]_q[b+1]_q}-\frac{1}{[a+1]_q}-\frac{1}{[b+1]_q}=q-1
    \end{align*}
    \begin{proof}[Proof of Lemma \ref{Lem2}]
    It follows that
    \begin{align*}
        &\frac{[a+b+2]_q}{[a+1]_q[b+1]_q}-\frac{1}{[a+1]_q}-\frac{1}{[b+1]_q}=\frac{[a+b+2]_q-[a+1]_q-[b+1]_q}{[a+1]_q[b+1]_q}\\
        &=\frac{q^{a+1}[b+1]_q-[b+1]_q}{[a+1]_q[b+1]_q}=\frac{q^{a+1}-1}{[a+1]_q}=q-1
    \end{align*}
    \end{proof}
\end{Lem}

Proposition \ref{Pro2} is proved by setting $a_j,b_{j^{\prime}}$ by $\lvert \alpha_j\rvert, \lvert \beta_{j^{\prime}} \rvert$. Therefore, $F_r(\lambda)$ holds for the column linear relations from Theorem \ref{linearrelation}, proving Step II.
\end{proof}
\textbf{Step III}. \textit{$F_r((n-s)^{s'})$ from the first step and the linear relations from second step determine $F_r(\lambda)$ for all $\lambda\subset \ell\times (n-s)$.}\\
This step was proved by Per \cite{alexandersson2021}[Proposition 28], hence we proved Theorem \ref{rec} (rectangular lemma).\\
\end{proof}

Let ${\lambda}=(\lambda_1,\cdots,\lambda_\ell)\subset \ell\times (n-s)$ with $\ell \leq s \leq(n-s)$. It is known that $X_{\lambda}$ is a linear combination of $\{e_n, e_{n-1,1}, \cdots, e_{n-\ell,\ell}\}$. In addition, $X_{((n-s)^i)}\in Span\{e_n, \cdots, e_{n-i,i}\}$ and when $i\neq j$, $X_{((n-s)^i)}$ and $X_{((n-s)^j)}$ are linearly independent because $e_{n-i,i}$ and $e_{n-j,j}$ are linearly independent. There are some explicit formulas of coefficients of $X_{((n-s)^i)}$ with $e$-basis in \cite{cho2019positivity,abrue2020modularlaw}. Without using the formulas from \cite{cho2019positivity,abrue2020modularlaw}, we found the other explicit formula for $X_{\lambda}$ with $e$-basis in combinatorial way which is similar to Theorem \ref{rec}. 

\begin{Thm}\label{e}
Let ${\lambda}=(\lambda_1,\cdots,\lambda_\ell)\subset \ell \times (n-\ell)$ for some $\ell\leq n/2$.
Then $X_\lambda$ is the same as \[\sum_{r=0}^{\ell} [n-2r]_q[n-r-\ell-1]_q![r]_q! G_r(\lambda) \cdot e_{n-r, r} \; ,\]
where
$G_r(\lambda)$'s are unimodal polynomials in $q$ with non-negative coefficients. Explicitly, we have
\[G_r(\lambda)=\sum_{I\in {[\ell] \choose r}}\tilde{q^{d_{\lambda,I}}}c_I \cdot f(\ell,n-r-1,{\lambda},I)\;,\]
where the each term appearing in the right-hand side is a non-negative polynomial with the same center. 
 $\tilde{q^{d_{\lambda, I}}}$ is defined by
    \[\tilde{q^{d_{\lambda, I}}}\coloneqq r(n-r)-\sum_{i=1}^{r}(\ell-r+i-a_i)-\sum_{i=1}^{r}\lambda_{\ell-a_i+1}\]
where $I=\{a_1, \cdots, a_r\}.$
\end{Thm}

Note that we do need $\ell \leq n-\ell$ since $[n-r-\ell-1]_q!$ is well-defined when $n-r-\ell-1$ is non-negative. However, it is known that $X_\lambda=X_{\lambda'}$ where $\lambda'$ is the conjugate of $\lambda$ so that Theorem \ref{e} implies $e$-positivity of $X_\lambda$ for the abelian case.\\

Before we prove Theorem \ref{e}, we show that $G_r(\lambda)$ is a non-negative polynomial. If $f(\ell,n-r-1,\lambda,I)$ is nonzero, then $r$ must be less than or equal to the length of $\lambda$ (say $k$) by Lemma \ref{powerlemma}. Since $\lambda$ is contained in a rectangle, we have $\lambda_1\leq n-r$. Then $f(\ell,n-r-1,\lambda,I)$ is positive if $\lambda_1<n-r$ by Lemma \ref{fpositive}, but when $\lambda_1=n-r$, we cannot apply Lemma \ref{fpositive}. However, in that case, $I$ should be $\{\ell-r+1,\ell-r+2,\cdots,\ell\}$ and in that case $f(\ell,n-r-1,\lambda,I)$ is positive.
 Therefore, $G_r(\lambda)$ is non-negative polynomial.\\

Now, we prove the Theorem \ref{e}.
\begin{proof}
We need three steps to prove the Theorem \ref{e}: 
\begin{enumerate}[Step I]
\item Show that the theorem holds when $\lambda=((n-a)^{a})$ for all $0\leq a\leq n-1$.
\item Show that $G_r(\lambda)$ satisfies the row and column linear relations in Theorem \ref{linearrelation}.
\item Prove that $G_r((n-a)^{a})$ from the first step and the linear relations from second step determine $G_r(\lambda)$ for abelian $\lambda$.
\end{enumerate}

\textbf{Step I} \textit{The Theorem \ref{e} holds when $\lambda=((n-a)^{a})$ for all $0\leq a\leq n-1$.}
\begin{proof}
First of all, we need to show $G_r((n-a)^a)=0$ for $a\neq r$. It is enough to show  $f(\ell,n-r-1,(n-a)^a,I)=0$ for all $I$ when $a\neq r$. Let $I=\{a_1,\cdots,a_r\}$ and $\lambda=(n-a)^a$.



Define a partition $\mu=(\overbrace{n-a,\cdots,n-a}^{a},\overbrace{0,\cdots,0}^{\ell-a})$ with $\mu_{a+1}=\cdots = \mu_{\ell}=0$. If $G_r((n-a)^a)$ is nonzero, then $I$ must be $\{\ell-a+1,\cdots, \ell-a+r\}$ by Lemma \ref{powerlemma}. It follows that $a\geq r$. If $a>r$, then $f(\ell,n-r-1,\lambda,I)$ is zero because the term from the edge $(\ell-a+r+1,\ell+1)$ is $(n-r-1)-(n-a)-(a-r-1)=0$.\\

If $a=r$, $f(\ell,n-r-1,\lambda,I)$ is not zero if only if $I=\{\ell-r+1,\ell-r+2,\cdots, \ell\}$. In this case, we have $\tilde{q^{d_{\lambda,I}}}=c_I=1$. Therefore, then $G_r((n-r)^r)=[n-r]_q\cdots [n-r-r+1]_q\cdot[n-2r-1]_q\cdots [n-r-\ell]_q$, which implies
\begin{align*}
&\sum_{r=0}^{[\frac{n}{2}]} [n-2r]_q[n-r-\ell-1]_q![r]_q! G_r((n-a)^a) \cdot e_{n-r, r}\\
=&[n-2a]_q[n-a-\ell-1]_q![a]_q! G_a((n-a)^a) \cdot e_{n-a, a}\\
    =&[n-2a]_q[n-a-\ell-1]_q![a]_q!\\
    &\times\big([n-a]_q\cdots [n-2a+1]_q\cdot[n-2a-1]_q\cdots [n-a-\ell]_q\big)e_{n-a,a}\\
    =&[n-a]_q![a]_q! e_{n-a,a}=X_{((n-a)^a)}.
\end{align*}

\end{proof}

We observe that the differences between $F_r(\lambda)$ and $G_r(\lambda)$ do not affect the proof of Step II and III, therefore the proof of Step II and III from the theorem \ref{e} is exactly the same with the proof of Step II and III from rectangular lemma (Theorem \ref{rec}). Thus the proof of Theorem \ref{e} is completed.
\end{proof}

\section{Relation with rook placements}

In this section, we describe the relationship between the previous section and rook placements.\\

Let $m_1, m_2$ be positive integers with $m_1\leq m_2$. A rook placement in a $m_1 \times m_2$ board is a set of $m_1$ rooks placed in cells of the board such that there is no column or row containing more than one rook. Let $\lambda$ be a partition contained in $m_1\times m_2$. Each rook placement has a $\lambda$-weight defined as in \cite{dworkin1998}. The weight is the number of cells $e$ in $m_1 \times m_2$ board such that 
\begin{enumerate}\item there is no rook in $c$,
 \item there is no rook to the left of $c$,
\item if $c$ is in $\lambda$, then the rook on the same column of $c$ is in $\lambda$ and below $c$,
\item if $c$ is not in $\lambda$, then the rook on the same column of $c$ is either in $\lambda$ or below $c$.
\end{enumerate}

Then \emph{$q$-hit number} $H^{m_1,m_2}_j(\lambda)$ is defined by $\sum q^{wt_\lambda(p)}$ where the sum is over all rook placements $p$ in $m_1\times m_2$ such that $j$ rooks are in $\lambda$ and $wt_\lambda$ is the $\lambda$-weight of $p$. If $m_1=m_2$, we denote $H^{m_1,m_2}_j(\lambda)$ by $H^{m_1}_j(\lambda)$.

Now, we are ready to describe the relationship between our work and $q$-hit numbers.

\begin{Thm} \label{qhit} Assume that $\lambda$ is contained in $\ell \times (n-s)$ where $\ell \leq n-s$. Then we have

$$X_\lambda(x,q)= { 1 \over [n-s]_q[n-s-1] \cdots [n-s-\ell+1]_q} \sum_{r=0}^\ell H^{\ell,n-s}_r(\lambda) X_{((n-s)^r)}(x,q).$$

\end{Thm}
\begin{rem}
As of writing this paper after studying the above relationship, authors found that Theorem \ref{qhit} is also noticed by Guay-Paquet (unpublished).
\end{rem}
\begin{proof}
The idea of the proof is the same as the proof of Theorem \ref{rec}: Prove that the right-hand side of the theorem satisfies Step I, II, and, III. Step I is straight-forward since $H^{\ell,n-s}_r((n-s)^a)$ is zero if $a \neq r$, and $[n-s]_q[n-s-1]_q\cdots [n-s-\ell-1]$ if $a=r$. Also, $H^{\ell,n-s}_r(\lambda)$ satisfies the row and column linear relations by \cite[Lemma 4.2]{abrue2020modularlaw}, which is Step II. Lastly, Step III is the same as Step III of Theorem \ref{rec}, which is \cite[Propsition 28]{alexandersson2021}. 
\end{proof}\label{F=H} 
By Theorem \ref{rec} and \ref{qhit}, we have the following:
\begin{Cor} \label{Fqhit} For $\lambda \subset \ell \times (n-\ell)$, we have
$$F_r^{\ell,n-s}(\lambda)= H^{\ell,n-s}_r(\lambda).$$

\end{Cor} 
Therefore, the Corollary \ref{Fqhit} and Theorem \ref{rec} provide a positive explicit formula for the $q$-hit number.

\subsection{Explicit formulas for $e$-positivity for abelian cases}
We start with the main theorem of Abreu-Nigro \cite{abrue2020modularlaw}.
\begin{Thm}[Abreu-Nigro, 2020] \label{an}
Assume $\lambda \subset \ell \times (n-\ell)$ for some $\ell$ and $\ell(\lambda)=k\leq \lambda_1$.
\begin{align*}X_\lambda(x,q)&=[k]_q! H^{n-k}_k(\lambda)\cdot e_{n-k,k}\\& + \sum_{r=0}^{k-1} q^r [r]_q![n-2r]_q H^{n-r-1}_r(\lambda)\cdot e_{n-r,r}.\end{align*}
\end{Thm}

Therefore, Corollary \ref{Fqhit} also provide a manifestly positive formula for $e$-expansion of $X_\lambda(x,q)$ when $\lambda$ is abelian. One can compare Theorem \ref{an} and Theorem \ref{e}, then the main corollary is essentially the same as Corollary \ref{Fqhit}. 

\section{Explicit formulas for $e$-positivity for some non-abelian cases}

In this section, we apply our rectangular lemma and Theorem \ref{e} to find a few explicit formulas for $e$-positivity of chromatic quasisymmetric functions when $\lambda$ is a partition contained in a rectangle except the first row of $\lambda$. In addition, we also suggest some conjectural formulas for $e$-positivity of chromatic quasisymmetric functions when $\lambda$ is not abelian.

\begin{figure}[h]
    \centering
	    \includegraphics[scale=0.5]{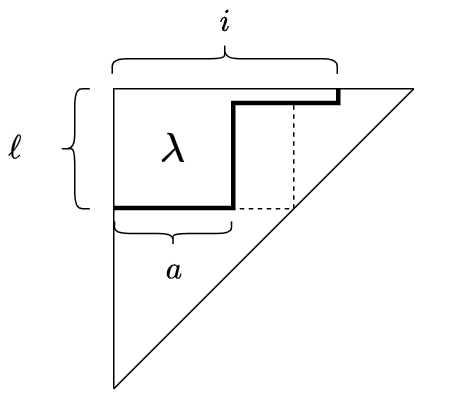}
	    \caption{${\lambda}=(i,a,\cdots,a)$}\label{nonabelian1}
\end{figure}

Let ${\lambda}=(i,a,\cdots,a)$ where $l({\lambda})=\ell$ and $a\leq n-\ell\leq i$. Then, $\lambda$ is a non-abelian partition but contained in a rectangle except the first row of $\lambda$, see figure \ref{nonabelian1}.

\begin{Thm} \label{nonlinear_n-1}
Let ${\lambda}=(i,a,\cdots,a)$ where $l({\lambda})=\ell$ and $a\leq n-\ell\leq i$. Then the coefficients of $e_{n-1,1}$ of $X_{\lambda}$ is the same as
\begin{align}
    &[n-3-a]_q\cdots[n-\ell-a]_q\cdot[n-\ell-2]_q!\nonumber\\
    &\text{\scalebox{0.94}{$\times \Big(q^{n-\ell-a}[n-1-i]_q[n-1]_q[a]_q[\ell-2]_q+q^{n-1-i}[i]_q[n-2-a]_q[n-\ell-1]_q\Big)$}}_{.}\label{4}
\end{align}

\end{Thm}
\begin{proof}
We will often use the linear relation from Theorem \ref{linearrelation2}.\\

Let ${\lambda}^{(i)}=(i,a,\cdots,a)$ be a length of $\ell$, where $n-\ell\leq i\leq n-1$. Then, by Theorem \ref{linearrelation2}, \[X_{\lambda^{(i)}}=\frac{[n-1-i]_q}{[\ell-1]_q} X_{\lambda^{(n-\ell)}}+\frac{q^{n-1-i}[i-(n-\ell)]_q}{[\ell-1]_q} X_{\lambda^{(n-1)}}.\]
Let $C_{\lambda^{(n-\ell)}}, C_{\lambda^{(n-1)}}$ be the coefficient of $e_{n-1,1}$ of $X_{\lambda^{(n-\ell)}}, X_{\lambda^{(n-1)}}$, respectively. By Theorem \ref{rec}, we have
\begin{align*}
    C_{\lambda^{(n-\ell)}}&=[n-2]_q[n-3-a]_q\cdots[n-\ell-a]_q\cdot[n-\ell-2]_q!\\
 &\times\Big(q^{n-\ell-a}[\ell]_q[a]_q[\ell-2]_q+q^{\ell-1}[n-\ell-a]_q[n-\ell-1]_q\Big)
\end{align*}
Also, $C_{\lambda^{(n-1)}}$ is equal to the coefficient of $e_{n-1}$ of $X_{((a)^{\ell-1})}$, hence we have   

 \[ C_{\lambda^{(n-1)}}=[n-1]_q[n-2-a]_q\cdots[n-\ell-a]_q \cdot [n-\ell-1]_q!,\]
Therefore, the coefficients of $e_{n-1,1}$ of $X_{\lambda^{(i)}}$ is the same as
\begin{align}
C_{\lambda^{(i)}}&=\frac{[n-3-a]_q\cdots[n-\ell-a]_q\cdot[n-\ell-2]_q!}{[\ell-1]_q}\nonumber\\
&\text{\scalebox{0.98}{$\times\bigg\{ [n-1-i]_q[n-2]_q\Big(q^{n-\ell-a}[\ell]_q[a]_q[\ell-2]_q+q^{\ell-1}[n-\ell-a]_q[n-\ell-1]_q\Big)$}}\nonumber\\
&\qquad+q^{n-1-i}[i-(n-\ell)]_q[n-1]_q[n-2-a]_q[n-\ell-1]_q\bigg\}\label{5}
\end{align}
We need to show the formula from the Theorem 4.1 is exactly the same as $C_{\lambda^{(i)}}$, therefore it is enough to show \[(\ref{4})-(\ref{5})=\frac{[n-3-a]_q\cdots[n-\ell-a]_q\cdot[n-\ell-2]_q!}{[\ell-1]_q}\times\mathcal{C}=0.\]
where $\mathcal{C}$ is as follows:

\begin{align*}
    \mathcal{C}&=\Big\{q^{n-\ell-a}[n-1-i]_q[a]_q[\ell-2]_q \big([\ell-1]_q[n-1]_q-[\ell]_q[n-2]_q \big)\\
    &\quad -q^{\ell-1}[n-\ell-a]_q[n-\ell-1]_q[n-2]_q[n-1-i]_q\\
    &\quad +q^{n-1-i}[n-2-a]_q[n-\ell-1]_q \big([i]_q[\ell-1]_q-[i-(n-\ell)]_q[n-1]_q \big)\Big\}\\
    &=\Big\{q^{n-\ell-a}[n-1-i]_q[a]_q[\ell-2]_q \cdot (-q^{\ell-1}[n-\ell-1]_q)\\
    &\quad-q^{\ell-1}[n-\ell-a]_q[n-\ell-1]_q[n-2]_q[n-1-i]_q\\
    &\quad+q^{n-1-i}[n-2-a]_q[n-\ell-1]_q\cdot q^{i-(n-\ell)}[n-1-i]_q[n-\ell]_q\Big\}\\
    &=[n-\ell-1]_q[n-1-i]_q\Big\{q^{n-a-1}[a]_q[\ell-2]_q\\
    &\quad -q^{\ell-1}([n-\ell-a]_q[n-2]_q-[n-2-a]_q[n-\ell]_q)\Big\}\\
    &=[n-\ell-1]_q[n-1-i]_q\Big(q^{n-a-1}[a]_q[\ell-2]_q-q^{\ell-1}(q^{n-\ell-a}[a]_q[\ell-2]_q)\Big)\\
    &=0
\end{align*}
Therefore, 
\begin{align*}
    C_{\lambda^{(i)}}&=\text{\scalebox{0.92}{$[n-3-a]_q\cdots[n-\ell-a]_q\cdot[n-\ell-2]_q!$}}\\
    &\times\Big(\text{\scalebox{0.92}{$q^{n-\ell-a}[n-1-i]_q[n-1]_q[a]_q[\ell-2]_q+q^{n-1-i}[i]_q[n-2-a]_q[n-\ell-1]_q$}}\Big)_{.}
\end{align*}
\end{proof}

\begin{Thm}
Let ${\lambda}=(i,a,\cdots,a)$ where $l({\lambda})=\ell$ and $a\leq n-\ell\leq i$. Then the coefficients of $e_{n-2,2}$ of $X_{\lambda}$ is the same as
\begin{align*}
    &[n-1-i]_q[n-4]_q[n-4-a]_q\cdots[n-\ell-a]_q\cdot[n-\ell-3]_q!\nonumber\\
    &\times \Big(q^{n-2-a}[2]_q[a]_q[n-\ell-a]_q[n-\ell-2]_q+q^{2(n-\ell-a)}[\ell]_q[a]_q[a-1]_q[\ell-3]_q\Big)_{.}
\end{align*}
\end{Thm}
\begin{proof}
We also use the linear relation from Theorem \ref{linearrelation2}.\\
Let ${\lambda}^{(i)}$ be $(i,a,\cdots,a)$ of length $\ell$, where $n-\ell\leq i\leq n-1$. Then, by Theorem \ref{linearrelation2}, we have \[X_{\lambda^{(i)}}=\frac{[n-1-i]_q}{[\ell-1]_q} X_{\lambda^{(n-\ell)}}+\frac{q^{n-1-i}[i-(n-\ell)]_q}{[\ell-1]_q} X_{\lambda^{(n-1)}}\]
Let $D_{\lambda^{(n-\ell)}}, D_{\lambda^{(n-1)}}$ be the coefficient of $e_{n-2,2}$ of $X_{\lambda^{(n-\ell)}}, X_{\lambda^{(n-1)}}$, respectively.
Since $D_{\lambda^{(n-1)}}=0$ and 
\begin{align*}
    D_{\lambda^{(n-\ell)}}=&[n-4]_q[n-\ell-3]_q![2]_q\\
    &\times\Big\{(q^{n-2-a}[\ell-1]_q[a]_q[n-\ell-a]_q[n-4-a]_q\cdots[n-\ell-a]_q[n-\ell-2]_q\\
    &\quad +q^{2(n-\ell-a)}\frac{[\ell]_q[\ell-1]_q}{[2]_q}[a]_q[a-1]_q[\ell-3]_q[n-4-a]_q\cdots[n-\ell-a]_q\Big\}_{,}
\end{align*}
Therefore, the coefficients of $e_{n-2,2}$ of $X_{\lambda^{(i)}}$ is the same as
\begin{align*}
    D_{\lambda^{(i)}}=&\frac{[n-1-i]_q}{[\ell-1]_q} \cdot D_{\lambda^{(0)}}\\
    =&\text{\scalebox{0.96}{$[n-1-i]_q[n-4]_q[n-\ell-3]_q!\cdot [n-4-a]_q\cdots[n-\ell-a]_q$}}\\
    &\text{\scalebox{0.94}{$\times \Big\{(q^{n-2-a}[2]_q[a]_q[n-\ell-a]_q[n-\ell-2]_q+q^{2(n-\ell-a)}[\ell]_q[a]_q[a-1]_q[\ell-3]_q\Big\}$}}_{.}
\end{align*}
\end{proof}

Now, we use the result from Theorem \ref{nonlinear_n-1} and rectangular lemma to find the coefficients of $e_{n-1,1}$ of $X_\lambda$ when $\lambda=(i, \overbrace{a, \cdots, a}^{\ell-1}, b_1, \cdots, b_p)$ where $b_p \leq \cdots \leq b_1 \leq a \leq n-\ell \leq i$.

\begin{Thm} \label{nonlinear1}
Let $\lambda=(i, \overbrace{a, \cdots, a}^{\ell-1}, b_1, \cdots, b_p)$ where $b_p \leq \cdots \leq b_1 \leq a \leq n-\ell-p$ and $n-\ell \leq i$. Then the coefficients of $e_{n-1,1}$ of $X_\lambda$ is the same as
\begin{align*}
    &[n-3-a]_q\cdots [n-\ell-a]_q\cdot [n-\ell-p-2]_q!\\
    &\times \Big (q^{n-1-i}[i]_q[n-2-a]_q[n-\ell-p-1]_q\prod_{i=1}^{i=p}[n-\ell-i-b_i]_q\\
    &\quad +[n-1-i]_q[n-1]_q\cdot q^{n-\ell-a}\sum_{I\in {[p+1] \choose 1}}q^{d_{\nu,I}}f(p+1,n-\ell-1,\nu, I)[\ell+p-1-k]_q\Big)_{,}
\end{align*}
where $\nu=(a,b_1,\cdots, b_p) \subset p \times a \text{ and } I=\{k\}$. 
\end{Thm}
\begin{proof}
We need three steps to prove the Theorem \ref{nonlinear1}: 
\begin{enumerate}[Step I]
\item Show that the theorem holds when $\lambda^{(j)}=(i,\overbrace{a,\cdots,a}^{j-1})$ for all $\ell \leq j \leq n$.
\item Show that the theorem satisfies the row and column linear relations in Theorem \ref{linearrelation} when the first $\ell$ part of $\lambda$ is fixed.
\item Prove that the first step and linear relations in Theorem \ref{linearrelation} determine the theorem for all $\ell$ and $p$.
\end{enumerate}

The idea is that the proof is very similar to the proof of Theorem \ref{rec} even though partitions appearing in this theorem is not abelian. The reason why we can apply the same argument is that the row and column linear relations can be applied when we fix the first $\ell$ part of $\lambda$, due to the inequality $a\leq n-\ell-p\leq n-\ell \leq i$. By using those linear relations, one can write $X_{\lambda}$ as a linear combination of $X_\rho$ where $\rho$ is one of $\lambda^{(j)}$=$(i,a,a,\cdots,a)$ for $\ell-1\leq j \leq p$. \\

 Also, note that the proof of Step III is essentially the same as the proof of Step III in Theorem \ref{rec} because $[\ell+p-1-k]_q$ is a quantity analogues to $c_I$. Thus, it is enough to show step I and step II.\\

\textbf{Step I} \textit{The Theorem 4.3 holds when $\lambda^{(j)}=(i,\overbrace{a,\cdots,a}^{j-1})$ for all $\ell \leq j \leq p$.}

\begin{proof}
Consider $\lambda^{(\ell+p')}$ for some $0\leq p'\leq p$. We have $b_1=\cdots b_{p'}=a$ and $b_{p'+1}=\cdots=b_p=0$. Then $f(p+1,n-\ell-1,\nu, I) \neq 0$ only when $I=\{k\}=\{p-p'+1\}$ by Lemma \ref{powerlemma} and we have
\begin{align*} f(p+1,n-\ell-1,\nu,I)&=[a]_q[n-\ell-1-a]_q\cdots[n-\ell-p'-a]_q\\
&\cdot [n-\ell-p'-2]_q[n-\ell-p'-3]_q \cdots [n-\ell-p-1]_q 
\end{align*}

and $d_{\nu,I}=1\cdot a - (p+1-k)- \nu_{p+2-k}=k-p-1=-p'$. Therefore
\begin{align*}
   &[n-3-a]_q\cdots [n-\ell-a]_q\cdot [n-\ell-p-2]_q!\\
 &\text{\scalebox{0.90}{$\times \Big (q^{n-1-i}[i]_q[n-2-a]_q [n-\ell-p-1]_q\prod_{i=1}^{p'}[n-\ell-i-a]_q\prod_{i=p'+1}^p[n-\ell-i]_q$}}\\
&\quad\text{\scalebox{0.90}{$+[n-1-i]_q[n-1]_q q^{n-\ell-a-p'}[a]_q\prod_{i=1}^{p'}[n-\ell-i-a]_q\prod_{i=p'+2}^{p+1}[n-\ell-i]_q$}}\Big)\\
&=[n-3-a]_q\cdots [n-\ell-p'-a]_q\cdot [n-\ell-p'-2]_q! [\ell+p-1-k]_q\\
 &\quad\text{\scalebox{0.90}{$\times \Big (q^{n-1-i}[i]_q[n-2-a]_q[n-\ell-p'-1]_q+q^{n-\ell-a-p'}[n-1-i]_q[n-1]_q [a]_q [\ell+p'-2]_q$}}\Big)\\
&=\text{the coefficients of} \; e_{n-1,1}\; \text{of} \; X_{\lambda^{(\ell+p')}} \text{by Theorem \ref{nonlinear_n-1}}.
\end{align*}
\end{proof}


\textbf{Step II} \textit{The Theorem satisfies the row and column linear relations in Theorem \ref{linearrelation}.}

\begin{proof}
First of all, instead of considering $(i, \overbrace{a, \cdots, a}^{\ell-1}, b_1, \cdots, b_p)$, we only need to consider $\nu=(a,b_1,\cdots, b_p)$ for the row and column relations. We divide the formula into two parts, $A$ and $B$, where $A$ is the term with $q^{n-1-i}[i]_q$ and $B$ is the terms with $f(p+1,n-\ell-1,\nu, I)$. We want to show each part satisfies the row and column relations separately.

For the row relation, let a partition $\nu$ with $j$ such that $b_{j}+2\leq b_{j-1}$, and $\nu^0=\nu$, $\nu^1, \nu^2$ be partitions defined by $\nu^a_k=\nu_k$ if $k\neq j$ and $\nu^a_j=\nu_j+e$ for $e=0,1,2$. For $A$ and $B$, since all other $q$-integers are exactly the same except the $q$-integers with $b_j$, it is enough to check the terms with $b_{j}$. For $A$, we have $[n-\ell-j-b_j]_q+q\cdot[n-\ell-j-b_j-2]_q=[2]_q[n-\ell-j-b_j-1]_q$ and for $B$, we have $q^{d_{\nu^0,I}}f(p+1,n-\ell-1,\nu^0,I)+q\cdot q^{d_{\nu^2,I}}f(p+1,n-\ell-1,\nu^2, I)=[2]_q q^{d_{\nu^1,I}}f(p+1,n-\ell-1,\nu^1,  I)$, which is exactly the same with the proof of Step II in Rectangular Lemma.

For the column relation, let a partition $\nu$ with $j$ such that $\nu_{j+1}=\nu_j \leq \nu_{i-1}-1$, and let $\nu^{(0)}=\nu, \nu^{(1)}, \nu^{(2)}$ be partitions defined by $\nu^{(1)}_k=\nu_k \text{ if } k\neq j$, $\nu^{(1)}_j=\nu_j+1$, and $\nu^{(2)}_k=\nu_k \text{ if } k \neq j, j+1$, $\nu^{(2)}_j=\nu_j+1=\nu^{(2)}_{j+1}$. For $A$ and $B$, since all other $q$-integers are exactly the same except the $q$-integers with $b_j$ and $b_{j+1}$, it is enough to check the terms with $b_j$ and $b_{j+1}$. For $A$, $[n-\ell-j-b_j]_q[n-\ell-j-1-b_{j+1}]_q+q\cdot[n-\ell-j-1-b_j]_q[n-\ell-j-2-b_{j+1}]_q=[n-\ell-j-1-b_j]_q\big([n-\ell-j-b_j]_q+q\cdot[n-\ell-j-2-b_{j+1}]_q\big)=[n-\ell-j-1-b_j]_q\cdot [2]_q[n-\ell-j-1-b_j]_q=[2]_q[n-\ell-j-1-b_j]_q[n-\ell-j-1-b_{j+1}]_q$. For $B$, we can apply the proof of Proposition \ref{Pro1} in rectangular lemma.
\end{proof}
Therefore, the proof is completed.
\end{proof}

\begin{figure}[h]
    \centering
	    \includegraphics[scale=0.65]{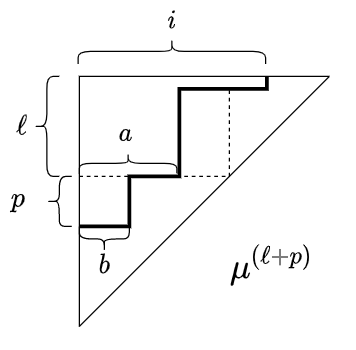}
	    \caption{${\mu^{(\ell+p)}}$} \label{nonabelian2}
\end{figure}

We can apply Theorem \ref{nonlinear1} to find the coefficient of $e_{n-1,1}$ of $X_\mu$ when $\mu=(i,\overbrace{a,\cdots,a}^\text{$\ell-1$},\overbrace{b,\cdots,b}^\text{$p$})$ with $0\leq b\leq a\leq n-\ell \leq i \leq n-1$. See figure \ref{nonabelian2}.

\begin{Cor}
Let ${\mu^{(\ell+p)}}=(i,\overbrace{a,\cdots,a}^\text{$\ell-1$},\overbrace{b,\cdots,b}^\text{$p$})$ with $0\leq b\leq a\leq n-\ell \leq i \leq n-1$.  Then the coefficients of $e_{n-1,1}$ of $X_{\mu^{(\ell+p)}}$ is
\begin{align*}
    &\text{\scalebox{0.85}{$[n-3-a]_q\cdots[n-\ell-a]_q\cdot[n-\ell-p-2]_q!$}}\\
    &\text{\scalebox{0.85}{$\times\Big(q^{n-\ell-a}[n-1-i]_q[n-1]_q[a-b]_q[n-\ell-2-b]_q\cdots[n-\ell-p-b]_q[n-\ell-p-1]_q[\ell-2]_q$}}\\
    &\quad\text{\scalebox{0.82}{$+q^{n-\ell-p-b}[n-1-i]_q[n-1]_q[b]_q[n-\ell-1-a]_q[n-\ell-2-b]_q\cdots[n-\ell-p-b]_q[\ell+p-2]_q$}}\\
    &\text{\scalebox{0.85}{$\quad+q^{n-1-i}[i]_q[n-2-a]_q[n-\ell-1-b]_q\cdots[n-\ell-p-b]_q[n-\ell-p-1]_q\Big)$}}_{.}
\end{align*}

\end{Cor}

\begin{proof}
We can apply Theorem \ref{nonlinear1} to prove the Corollary. Since $b_i=b$ for $\forall i\in [p]$, $f(p+1,n-\ell-1,\nu, I)$ only survived when $I=\{1\}$ or $I=\{p+1\}$, where $\nu=(a, \overbrace{b, \cdots, b}^{p})$. The proof is completed.
\end{proof}

We still try to find the formula for the coefficients of $e_{n-2,2}$ of $X_{\lambda^{(\ell+p)}}$ where $\lambda^{(\ell+p)}=(i, \overbrace{a, \cdots, a}^{\ell-1}, b_1, \cdots, b_p)$, $b_p \leq \cdots \leq b_1 \leq a \leq n-\ell \leq i$, $p\leq n-\ell$, and $(b_1, \cdots, b_p)$ fits in $p \times (n-\ell-p)$. Since it has too many terms, we first consider to find the coefficients of $e_{n-2,2}$ of $X_{\lambda}$ when $\lambda=(i,\overbrace{a,\cdots,a}^\text{$\ell-1$},\overbrace{b,\cdots,b}^\text{$p$})$ as follows.

\begin{Con} \label{31}
Let ${\lambda}=(i,\overbrace{a,\cdots,a}^\text{$\ell-1$},\overbrace{b,\cdots,b}^\text{$p$})$ with $0\leq b\leq a\leq n-\ell \leq i \leq n-1$.

Then the coefficients of $e_{n-2,2}$ of $X_{\lambda}$ when $p\geq 2$ is the same as
\begin{align*}
    &\text{\scalebox{0.78}{$[n-1-i]_q[n-4][n-4-a]_q\cdots[n-\ell-a]_q\cdot[n-\ell-3-b]_q\cdots[n-\ell-p-b]_q\cdot[n-\ell-p-3]_q!$}}\\
    &\times \Big(\text{\scalebox{0.85}{$q^{n-2-a}[2]_q[a-b]_q[n-\ell-a]_q[n-\ell-2-b]_q[n-\ell-p-1]_q[n-\ell-p-2]_q$}}\\
    &+\text{\scalebox{0.8}{$q^{n-2-b}[2]_q[b]_q[n-\ell-1-a]_q[n-\ell-2-b]_q[n-\ell-p-a]_q[n-\ell-p-2]_q$}}\\
    & +\text{\scalebox{0.8}{$q^{2(n-\ell-a)}[a-b]_q[a-b-1]_q[n-\ell-p-1]_q[n-\ell-p-2]_q[\ell]_q[\ell-3]_q$}}\\
    &+\text{\scalebox{0.78}{$q^{2n-2\ell-a-b-p-1}[2]_q[b]_q[a-b]_q[n-\ell-1-a]_q[n-\ell-p-2]_q\big([\ell+p-1]_q[\ell-2]_q+q^{\ell-1}[p-2]_q\big)$}}\\
    &+\text{\scalebox{0.8}{$q^{2(n-\ell-p-b)}[b]_q[b-1]_q[n-\ell-1-a]_q[n-\ell-2-a]_q[\ell+p]_q[\ell+p-3]_q$}}\Big)_{.}
\end{align*}
\end{Con}
Note that when $p=2$, there are no terms with $[n-\ell-3-b]_q\cdots[n-\ell-p-b]_q$ from the first row of Conjecture \ref{31}.

When $p=1$, the coefficients of $e_{n-2,2}$ of $X_{\lambda}$ is the same as
\begin{align*}
    &[n-1-i]_q[n-4][n-4-a]_q\cdots[n-\ell-a]_q\cdot[n-\ell-4]_q!\\
    &\times \Big(q^{n-2-a}[2]_q[a-b]_q[n-\ell-a]_q[n-\ell-2]_q[n-\ell-3]_q\\
    &+q^{n-2-b}[2]_q[b]_q[n-\ell-1-a]_q[n-\ell-1-a]_q[n-\ell-3]_q\\
    &+q^{2(n-\ell-a)}[a-b]_q[a-1]_q[n-\ell-3]_q[\ell]_q[\ell-3]_q\\
    &+q^{2n-2\ell-a-b-2}[b]_q[a-1]_q[n-\ell-1-a]_q[\ell+1]_q[\ell-2]_q\Big)_{.}
\end{align*}

\bibliographystyle{amsplain}
\bibliography{sk1}

\providecommand{\bysame}{\leavevmode\hbox to3em{\hrulefill}\thinspace}
\providecommand{\MR}{\relax\ifhmode\unskip\space\fi MR }
\providecommand{\MRhref}[2]{%
  \href{http://www.ams.org/mathscinet-getitem?mr=#1}{#2}
}
\providecommand{\href}[2]{#2}
\begin{thebibliography}{10}

\bibitem{abrue2020modularlaw}
Antonio~Nigro Alex~Abreu, \emph{Chromatic symmetric functions from the modular
  law}, arXiv:2006.00657.

\bibitem{alexandersson2021}
Per Alexandersson, \emph{$llt$ polynomials, elementary symmetric functions and
  melting lollipops}, Journal of Algebraic Combinatorics \textbf{52} (2021),
  299--325.

\bibitem{alexandersson2018llt}
Per Alexandersson and Greta Panova, \emph{${LLT}$ polynomials, chromatic
  quasisymmetric functions and graphs with cycles}, Discrete Mathematics
  \textbf{341} (2018), no.~12, 3453--3482.

\bibitem{BW1989}
Anders Björner and Michelle L.~Wachs, \emph{$q$-hook length formulas for
  forests}, Journal of Combinatorial Theory Series A \textbf{52} (1989), no.~2,
  165--187.

\bibitem{carlsson2018proof}
Erik Carlsson and Anton Mellit, \emph{A proof of the shuffle conjecture},
  Journal of the American Mathematical Society \textbf{31} (2018), no.~3,
  661--697.

\bibitem{cho2020bounce3}
Soojin Cho and Jaehyun Hong, \emph{Positivity of chromatic symmetric functions
  associated with hessenberg functions of bounce number 3}, arXiv:1910.07308.

\bibitem{cho2019positivity}
Soojin Cho and JiSun Huh, \emph{On e-positivity and e-unimodality of chromatic
  quasi-symmetric functions}, SIAM Journal on Discrete Mathematics \textbf{33}
  (2019), no.~4, 2286--2315.

\bibitem{dworkin1998}
M.~Dworkin, \emph{An interpretation for garsia and remmel’s q-hit numbers},
  J. Combin. Theory Ser. A \textbf{81} (1998), no.~2, 149--175.

\bibitem{gasharov1996incomparability}
Vesselin Gasharov, \emph{Incomparability graphs of (3+ 1)-free posets are
  s-positive}, Discrete Mathematics \textbf{157} (1996), no.~1-3, 193--197.

\bibitem{guay2013modular}
Mathieu Guay-Paquet, \emph{A modular relation for the chromatic symmetric
  functions of (3+ 1)-free posets}, arXiv preprint arXiv:1306.2400 (2013).

\bibitem{haglund2005combinatorial}
James Haglund, Mark Haiman, Nicholas Loehr, Jeffrey~B Remmel, Alexander
  Ulyanov, et~al., \emph{A combinatorial formula for the character of the
  diagonal coinvariants}, Duke Mathematical Journal \textbf{126} (2005), no.~2,
  195--232.

\bibitem{harada2019cohomology}
Megumi Harada and Martha~E Precup, \emph{The cohomology of abelian hessenberg
  varieties and the stanley--stembridge conjecture}, Algebraic Combinatorics
  \textbf{2} (2019), no.~6, 1059--1108.

\bibitem{huh2020melting}
JiSun Huh, Sun-Young Nam, and Meesue Yoo, \emph{Melting lollipop chromatic
  quasisymmetric functions and schur expansion of unicellular ${LLT}$
  polynomials}, Discrete Mathematics \textbf{343} (2020), no.~3, 111728.

\bibitem{Jangsoo2014}
Jang~Soo Kim, Karola Mészáros, Panova Greta, and David B.~Wilson, \emph{Dyck
  tilings, increasing trees, descents, and inversions}, Journal of
  Combinatorial Theory, Series A \textbf{122} (2014), 9--27.

\bibitem{lascoux1997ribbon}
Alain Lascoux, Bernard Leclerc, and Jean-Yves Thibon, \emph{Ribbon tableaux,
  hall--littlewood functions, quantum affine algebras, and unipotent
  varieties}, Journal of Mathematical Physics \textbf{38} (1997), no.~2,
  1041--1068.

\bibitem{lee2018linear}
Seung~Jin Lee, \emph{Linear relations on ${LLT}$ polynomials and their k-schur
  positivity for k= 2}, Journal of Algebraic Combinatorics \textbf{53} (2021),
  973--990.

\bibitem{shareshian2016chromatic}
John Shareshian and Michelle~L Wachs, \emph{Chromatic quasisymmetric
  functions}, Advances in Mathematics \textbf{295} (2016), 497--551.

\bibitem{Shi2012}
Keiichi Shigechi and Paul Zinn-Justin, \emph{Path representation of maximal
  parabolic kazhdan–lusztig polynomials}, Journal of Pure and Applied Algebra
  \textbf{216} (2012), no.~11, 2533--2548.

\bibitem{stanley1995symmetric}
Richard~P Stanley, \emph{A symmetric function generalization of the chromatic
  polynomial of a graph}, Advances in Mathematics \textbf{111} (1995), no.~1,
  166--194.

\bibitem{stanley1993immanants}
Richard~P Stanley and John~R Stembridge, \emph{On immanants of jacobi-trudi
  matrices and permutations with restricted position}, Journal of Combinatorial
  Theory Series A \textbf{62(2)} (1993), 261--279.

\bibitem{sun2015polynomials}
Hua Sun, Yi~Wang, and Hai~Xia Zhang, \emph{Polynomials with palindromic and
  unimodal coefficients}, Acta Mathematica Sinica, English Series \textbf{31}
  (2015), no.~4, 565--575.

\end{thebibliography}

\end{document}